\documentclass[reqno,twoside]{amsart}
\usepackage{amsmath, amssymb}

\theoremstyle{plain}
\newtheorem{theorem}{Theorem}
\newtheorem{lemma}[theorem]{Lemma}
\newtheorem{proposition}[theorem]{Proposition}
\newtheorem{corollary}[theorem]{Corollary}
\newtheorem{conjecture}[theorem]{Conjecture}

\theoremstyle{definition}
\newtheorem{definition}[theorem]{Definition}
\newtheorem{example}[theorem]{Example}
\newtheorem{condition}[theorem]{Condition}
\newtheorem{examples}[theorem]{Examples}

\theoremstyle{remark}
\newtheorem{remark}[theorem]{Remark}

\numberwithin{equation}{section} 
\numberwithin{theorem}{section}
\newcommand{\thref}[1]{Theorem~{\rm\ref{#1}}}

\newcommand{\leref}[1]{Lemma~{\rm\ref{#1}}}
\newcommand{\coref}[1]{Corollary~{\rm\ref{#1}}}
\newcommand{\condref}[1]{Condition~{\rm\ref{#1}}}
\newcommand{\deref}[1]{Definition~{\rm\ref{#1}}}
\newcommand{\exref}[1]{Example~{\rm\ref{#1}}}
\newcommand{\reref}[1]{Remark~{\rm\ref{#1}}}
\newcommand{\seref}[1]{Section~{\rm\ref{#1}}}
\newcommand\lbb[1]{\label{#1}}
\newcommand{\CC}{\mathbb{C}} 
\newcommand{\ZZ}{\mathbb{Z}}
\DeclareMathOperator{\fil}{F} 
\def\M{\mathcal{M}} 
\def\gg{\mathfrak{g}}
\def\bb{\mathfrak{b}} 
\def\hh{\mathfrak{h}} 
\newcommand{\WW}{\mathfrak{W}}
\newcommand{\eps}{\varepsilon}
\renewcommand\geq{\geqslant}
\renewcommand\leq{\leqslant} 
\DeclareMathOperator{\rk}{rk}
\DeclareMathOperator{\codim}{codim}
\DeclareMathOperator{\ann}{Ann}
\DeclareMathOperator{\Span}{Span}
\DeclareMathOperator{\gr}{gr}
\DeclareMathOperator{\ad}{ad}
\DeclareMathOperator{\id}{id}
\DeclareMathOperator{\Lin}{Lin}
\DeclareMathOperator{\Hom}{Hom}
\DeclareMathOperator{\End}{End}
\DeclareMathOperator{\Der}{Der} 
\DeclareMathOperator{\gk}{GKdim}
\DeclareMathOperator{\cur}{Cur}
\DeclareMathOperator{\dif}{Diff}
\DeclareMathOperator{\coef}{\mathcal{A}}
\DeclareMathOperator{\cend}{Cend}
\DeclareMathOperator{\chom}{Chom}
\renewcommand\Im{{\mathrm{Im}}\,}
\def\ot{\otimes}
\def\textsum{{\textstyle\sum}}
\def\smash{\,\sharp\,}
\def\ti#1{\widetilde{#1}}
\def\wti{\widetilde}
\def\De{\Delta}
\def\d{\partial}

\begin{document}

\title[Associative Pseudoalgebras]{Unital Associative Pseudoalgebras and 
Their Representations}

\author[A.~Retakh]{Alexander Retakh}

\address{Department of Mathematics
\\Yale University
\\New Haven, CT 06520}
\email{retakh@math.yale.edu}

\date\today

\begin{abstract}
Pseudoalgebras, introduced in \cite{BDK}, are multi-dimensional analogues
of conformal algebras, which provide an axiomatic description of the
singular part of the operator product expansion.

Our main interest in this paper is the pseudoalgebra $\cend_n$, which is 
the analogue of an algebra of endomorphisms of a finite module.  We study
its algebraic properties.  In particular, we introduce the class of 
unital pseudoalgebras and describe their structure and representations.  
Also, we classify pseudoalgebras algebraically similar to $\cend_n$.
\end{abstract}

\maketitle


\section*{Introduction}\lbb{sec.int}

Recent years saw several new approaches in the theory of vertex algebras.
One was to provide an axiomatic description of the ``singular'' part of
the vertex algebra (which came to be called a conformal algebra), thus
describing the vertex algebra as a highest-weight module (\cite{K1}).
The other strives for a coordinate-less description of vertex algebras, in
particular leading to the construction of analogous structures over any
complex curve (\cite{BD}, see also \cite{G} for an exposition and 
\cite{HL, F} for the relation to vertex algebras).  The latter approach is 
based on representation of algebras in pseudotensor categories.  For 
example, conformal algebras are algebras from the category of left
$\CC[\d]$-modules.

Using the language of pseudotensor categories, a natural generalization of
conformal algebras was introduced in \cite{BDK}.  These objects, called
pseudoalgebras, are also related to the differential Lie algebras of Ritt
and Hamiltonian formalism in the theory of nonlinear evolution equations.
In \cite{BDK} a full classification of semisimple finite Lie
pseudoalgebras was obtained, together with a foundation for their
representation theory.  One of the main objects there is the pseudoalgebra
of pseudolinear operators of a finite module.  Its study was the main
inspiration for this paper.

Our ground field is $\CC$; the reader may replace it with the favorite
algebraically closed field of zero characteristic.

\subsection{Pseudoalgebras} A pseudoalgebra is a left module $R$ over a
Hopf algebra $H$ together with an ``operation'' $R\ot R\to (H\ot H)\ot_H
R$ denoted $a*b$.  Here the tensor product over $H$ is defined via the
comultiplication $\De: H\to H\ot H$.

The Hopf structure allows to define the $x$-products $R\ot R\to R$ denoted
$a_xb$ for every $x$ in the dual algebra $H^*=\Hom(H,\CC)$. We have
\begin{equation}\lbb{eq.pseudomult}
a*b=\textsum_i (S(h_i)\ot 1)\ot_H(a_{x_i}b),\quad\text{for dual bases
$\{h_i\}$, $\{x_i\}$ of $H,H^*$}.
\end{equation}
Thus, a pseudoalgebra can be defined as an $H$-module with a collection of
bilinear operations indexed by $H^*$.

Consider the annihilation algebra of a pseudoalgebra $R$, 
$\coef(R)=H^*\ot_H R$. $\coef(R)$ possesses a natural left $H$-action (on
the first component of the tensor product).  For $a,b\in R$, 
$a*b=\textsum_i (f_i\ot g_i)\ot_H c_i$, multiplication is defined as 
\begin{equation}\lbb{eq.annihmult}
(x\ot_H a)(y\ot_H b)=\textsum_i (xf_i)(yg_i)\ot_H c_i.
\end{equation}
When $R$ is $H$-torsion free, the annihilation algebra $\coef(R)$ almost
completely describes $R$.  In particular, $R$ is a left $\coef(R)$-module:
$(x\ot_H a)\cdot b=a_xb$.

One can also define varieties of pseudoalgebras.  In the Lie and
associative cases, this requires extending the pseudoproducts to maps
$(H^{\ot 2}\ot_H R)\ot R\to H^{\ot 3}\ot_H R$ and $R\ot (H^{\ot 2}\ot_H
R)\to H^{\ot 3}\ot_H R$.  This is done via the pseudotensor structure on
the category of left $H$-modules.  In particular, for $a,b\in R$,
$a*b=\sum_i \beta_i\ot_H c_i$, we have
\begin{align*}
(\gamma\ot_H a)*b&=\textsum_i(\gamma\ot 1)(\De\ot\id)(\beta_i)\ot_H c_i,\\
a*(\gamma\ot_H b)&=\textsum_i(1\ot\gamma)(\id\ot\De)(\beta_i)\ot_H c_i.
\end{align*}
Then associativity is given by a familiarly looking equality in $H^{\ot
3}\ot_H R$:
\begin{equation*}
a*(b*c)=(a*b)*c.
\end{equation*}

Pseudomodules over pseudoalgebras of a particular variety are defined in a
similar way.

One of the most important examples of pseudoalgebras is $\cend_n$, a
pseudoalgebra of pseudolinear operators of a finite $H$-modules.  As the
category of $H$-modules plays a role of the category vector spaces in the
ordinary (non-pseudo) representation theory, $\cend_n$ is an analogue of
the algebra $\End_n$.  In particular, endowing a free $H$-module on $n$
generators with a structure of an $R$-pseudomodule for an associative 
pseudoalgebra $R$ is the same as providing a map $R\to\cend_n$.  As every
associative pseudoalgebra can be made into a Lie one, a similar statement
is true for Lie pseudoalgebras as well.

\subsection{Unital algebras and their representations} Any study of
pseudoalgebras is ultimately a study of corresponding annihilation
algebras.  The standard trick in the study of ordinary associative
algebras is to adjoin the identity; however, it is unclear if such an
operation can be performed on the pseudoalgebra level, i.e. if we will
still remain in the class of annihilation algebras.  Thus, it is necessary
to introduce some concept of ``identity'' for the pseudoalgebras
themselves.

Identity in ordinary algebras comes from an embedding of the ground field
$\CC$.  Similarly, for pseudoalgebras unitality is formulated as the
existence of an embedding of the finite $H$-pseudoalgebra that is a free
$H$-module of $\rk 1$ such that its action on the whole pseudoalgebra is
non-zero.  Such pseudoalgebras are called {\em unital}.  Of course, a
pseudoalgebra of all pseudolinear endomorphisms of some $H$-module (in
particular, $\cend_n$) is unital. The generator of the resulting
subalgebra produces a left identity when passing to the annihilation
algebra.

The structure of $\cend_n$ can be described explicitly.  Namely,
$\cend_n=H\ot H\ot\End_n(\CC)$ and the multiplication is completely
determined by that in the ordinary algebra $H\ot\End_n(\CC)$ and the
$H^*$-action on it.  More generally, let $A$ be an
$(H^*)^{cop}$-differential algebra (when $H$ is commutative, this simply
means that elements of $H^*$ acts on $A$ as differential operators, see
(\ref{eq.diffalgebra})) such that the action of the augmentation ideal of
$H^*$ is locally nilpotent. Then $A$ gives rise to a pseudoalgebra $\dif
A=H\ot A$ where the multiplication is defined via (\ref{eq.pseudomult})
with $(1\ot a)_{x_i}(1\ot b)=1\ot(a\cdot x_i(b))$ for $a,b\in A$. Such
pseudoalgebras are called {\em differential} (see \exref{ex.diff}). They 
can be easily made unital (by adjoining identity to $A$). In most cases 
the converse is also true.

\begin{theorem}\lbb{mainthunital}  A semisimple unital associative
pseudoalgebra is differential.
\end{theorem}

Passing from $\dif A$ to $A$ (whose structure is easier to understand than
that of the annihilation algebra) allows for the study and, in some sense,
complete classification of representations of $\dif A$. 

\begin{theorem}\lbb{mainthrep}  Let $V$ be a representation of
unital differential pseudoalgebra $R=\dif A$.  Then
$V=V^0\oplus V^1$, where $R*V^0=0$ and $V^1$ is constructed from a unitary
$A$-module.  Moreover, $V^1$ is irreducible (indecomposable) if and only
if $A$ is irreducible (indecomposable).
\end{theorem}

In particular,  the above theorem allows to classify representations of 
$\cend_n$.

\subsection{Classification results} Let us now try to answer the question 
what algebras are ``similar'' to $\cend_n$ (i.e. have a similar 
algebraic structure). Clearly, the necessary conditions must be
simplicity and unitality.  A finitness condition must be a requirement as
well; even though, $\cend_n$ is not finite over $H$. In \cite{Re1} we
classified conformal algebras of linear growth; however, this is not a
good condition for the case of general $H$.

Given a pseudoalgebra $R$ over a Hopf subalgebra $H'$ of $H$, it is easy
to lift it to an $H$-pseudoalgebra together with its modules; this
construction is called the current extension $\cur^{H'}_H R$ (see
\exref{ex.cur}).  Therefore, it is only natural to include in any
classification current extensions of $H'$-pseudoalgebras with the same
properties.  On the other hand, the growth of such current extensions does
not change, so it is possible to produce a current $H$-pseudoalgebra that is
``small'' compared to $H$ but has properties of the original algebra.  
Hence, we ought to rely not on comparison of growth but on some absolute
criterion which does not change with passing to the current extension.  
Since in the ordinary ($H=\CC$) case the algebras $\cend_n=\End_n(\CC)$
are finite, we shall require a presence of a certain finite current 
subalgebra.

As we remarked above, $\cend_n$ is an example of a differential
pseudoalgebra.  The multiplication in $\cend_n$ is characterized by the
action of $H^*$ on $H\ot \End_n(\CC)$ by differential operators.  This
action survives passing to the universal enveloping algebra of a central
extension of $\gg$.  In particular, let $0\to\CC{\tt c}\to\hat\gg\to\gg\to
1$ be determined by a cocycle $\phi$. $H^*$ naturally acts on
$U(\hat\gg)/(1-{\tt c})$ and, consequently, on its tensor product with
$\End_n(\CC)$ (for details see \exref{ex.weyl}).  This defines the
differential pseudoalgebra $\cend_n^\phi$.

\begin{theorem}\lbb{mainthclass}  Let $R$ be a simple unital associative
$H$-pseudoalgebra. Assume that its maximal unital current subalgebra 
having the  same pseudoidentity as $R$ is  finite as an $H$-module and 
simple.  Then  $R$ is either of
\begin{itemize}
\item $\cur\End_n(\CC)$, $n\geq 0$;
\item $\cur^{H'}_H\cend_n$, $H'$ a Hopf subalgebra of $H$, $n>0$;
\item $\cur^{H'}_H\cend^\phi_n$, $H'=U(\hh)$ a Hopf subalgebra of
$H$, $\phi\in H^2(\hh,\CC)$, $n>0$.
\end{itemize}
\end{theorem}

Remark that $R$ above is automatically finitely generated.

Even for a commutative $H$, $\cend_n$ and $\cend_n^\phi$ (which is
constructed from the Weyl algebra) have very different representations.
However, their algebraic properties are remarkably similar.  This
discrepancy, of course, does not occur in the ordinary ($H=\CC$) and the 
conformal ($H=\CC[\d]$) cases.

\begin{remark}\lbb{homolclass} I am unaware of another algebraic
classification where objects in a given class are parametrized by
$\bigcup_{\hh\subset\gg} H^2(\hh,\CC)$ for a finite-dimensional Lie
algebra $\gg$.  For an arbitrary finite group $G$, elements of
$\bigcup_{H\leq G} H^2(H,\CC)$ correspond to indecomposable modular
categories over $\text{Rep } G$ \cite{O} but this result can not be 
carried over to the case of $\text{Rep }\gg$.
\end{remark}

\subsection{Organization of the paper} We start by defining in
\seref{sec.prelim} the main objects of our study: pseudoalgebras and their
representations.  We also recall some useful facts about Hopf algebras.
This Section mainly follows a similar discussion in the first chapters of
\cite{BDK}.

\seref{sec.ex} is devoted to examples of pseudoalgebras.  In particular,
we will define the differential pseudoalgebras in \exref{ex.diff} and
discuss in greater detail the structure of $\cend_n$.

In \seref{sec.unital} we define unital pseudoalgebras and show in
\thref{unitisdiff} that under a certain technical condition such algebras
are differential over unital algebras.  \thref{mainthunital} follows
(\coref{semisimdiff}).  The proof itself expands that of Proposition 3.5
in \cite{Re1}.  We also briefly touch upon the classification of unital
algebras over cocommutative Hopf algebras.

In \seref{sec.rep} we describe representations of unital differential
pseudoalgebras.  In particular, we explicitly show how to construct a
unitary pseudomodule of a differential pseudoalgebra $\dif A$ from a
unitary $A$-module (\leref{tilderep}).  \thref{mainthrep} follows.  As a
corollary we derive the classification of indecomposable modules of the
conformal algebra $\cend_n$.  It was originally stated in \cite{K2};
however, our methods are different and do not refer to the Lie case.

Finally, in \seref{sec.simple} we prove \thref{mainthclass}.  This
essentially boils down to classifying certain filtered algebras with a
locally nilpotent $X^{cop}$-action.  When $X$ is cocommutative, these
algebras are differentiably simple; however, our methods are different
from other studies of differentiably simple algebras, as there are no
minimality conditions (see \cite{Bl} and references therein, and also 
\cite{Kh}). We also discuss in greater detail the relation between 
properties of pseudoalgebras and their growth and conjecture a stronger 
version of \thref{mainthclass}.

\subsection*{Acknowledgments} I am grateful to J.~T.~Stafford and
E.~Zelmanov for fruitful discussions and to V.~G.~Kac for providing an
updated version of \cite{K2}.


\section{A short survey of pseudoalgebras}\lbb{sec.prelim}

\subsection{Pseudotensor categories} Theory of pseudotensor categories was
developed in \cite{BD} as a way of expressing such notions as Lie
algebras, representations etc. in purely categorical terms.  The ultimate
goal is to define these notions for categories of modules that have an
interesting action on tensor products, e.g. $\mathcal{D}$-modules (as in
\cite{BD}) or modules over a Hopf algebra (as in \cite{BDK} or this
paper).

More details on pseudo-tensor categories can be found in
\cite[Chapter 1]{BD} and \cite[Chapter 4]{BDK}, we will get by with a
short presentation of main definitions and several examples.

Denote by $\mathcal{S}$ the category of finite non-empty sets with
surjective maps.  For a morphism $\pi:J\twoheadrightarrow I$ and $i\in I$,
we denote $\pi^{-1}(i)=J_i$.

\begin{definition}\lbb{pseudotensor} A {\em pseudotensor category} is a
class of objects $\M$ together with the following data:
\begin{itemize}
\item For any $I\in\mathcal{S}$, a family of objects $\{L_i\}_{i\in I}$
and an object $M$, one has the set of {\em polylinear maps}
$\Lin_I(\{L_i\}, M)$; the symmetric group $S_I$ acts on $\Lin_I(\{L_i\}$.
\item For any morphism $\pi:J\twoheadrightarrow I$ in $\mathcal{S}$,
the families of objects $\{L_i\}_{i\in I}$ and $\{N_j\}_{j\in J}$ and an
object $M$, there exists the composition map
\begin{equation*}
\Lin_I(\{L_i\}, M)\ot \bigotimes_{i\in I}\Lin_{J_i}(\{N_j\},
L_i)\to\Lin_J(\{N_j\}, M)\ni
\phi(\{\psi_i\}_{i\in I}).
\end{equation*}
\end{itemize}
This data satisfies the following properties:
\begin{description}
\item[Associativity] For a surjective map $K\twoheadrightarrow J$ and a
$K$-family of objects $\{P_k\}_{k\in K}$, 
one has
$\phi(\{\psi_i(\{\chi_j\})\}) =  
(\phi(\{\psi_i\}))(\{\chi_j\})\in\Lin_K(\{P_k\}, M)$,
given $\chi_j\in\Lin_{K_j}(\{P_k\}, N_j)$;
\item[Unit] For any object $M$, there exists an element
$\id_M\in\Lin(\{M\}, M)$ such that for any $\phi\in\Lin_I(\{L_i\}, M)$, 
one has $\id_M(\phi)=\phi(\{\id_{L_i}\})=\phi$;
\item[Equivariance] The compositions of polylinear maps are equivariant
with respect to the natural action of the symmetric group.
\end{description}
\end{definition}

\begin{examples}\lbb{ex.pseudotensor} 1. For the category $\mathcal{V}ec$
of vector spaces, put $\Lin_I(\{L_i\}, M)=\Hom(\ot_i L_i, M)$. The
symmetric group acts on $\Lin_I(\{L_i\}, M)$ by permuting the factors of
$\ot_i L_i$.

2. Let $H$ be a cocommutative bialgebra with a comultiplication $\De:H\to
H^{\ot 2}$ and $\M^l(H)$ its category of left modules.  This is a
symmetric tensor category; hence, it can be made into a pseudotensor
category: $\Lin_I(\{L_i\}, M)=\Hom_H(\ot_i L_i, M)$.

3.  We will introduce another pseudotensor structure on $\M^l(H)$.

Recall that $\De$ gives rise to a functor $\M^l(H)\to\M^l(H^{\ot 2})$,
$M\mapsto H^{\ot 2}\ot_H M$, where $H$ acts on $H^{\ot 2}$ via $\De$. This
may be generalized as follows. For every surjection $\pi:
J\twoheadrightarrow I$, define a functor $\De^{(\pi)}: \M^l(H^{\ot I})\to
\M^l(H^{\ot J})$, $M\mapsto H^{\ot J}\ot_{H^{\ot I}} M$, where $H^{\ot I}$
acts on $H^{\ot J}$ via the iterated comultiplication determined by $\pi$
(the $i$-th copy of $H$ is mapped into $H^{\ot J_i}$).  This is
well-defined because of coassociativity.

Denote the tensor product functor $\M^l(H)^I\to \M^l(H^{\ot I})$ by
$\boxtimes_{i\in I}$.  Then we can define a pseudotensor
category $\M^*(H)$ that has the same objects as $\M^l(H)$ but with
\begin{equation}\lbb{pseudotensorhopf}
\Lin_I(\{L_i\}, M)=\Hom_{H^{\ot I}} (\boxtimes_{i\in I}L_i, H^{\ot
I}\ot_H M).
\end{equation}
For $\pi:J\twoheadrightarrow I$, the composition of polylinear maps is
defined as follows:
\begin{equation}\lbb{pseudotensorhopfcomp}
\phi(\{\psi_i\})=\De^{(\pi)}(\phi)\circ(\boxtimes_i\psi_i).
\end{equation}

The symmetric group acts on $\Lin_I(\{L_i\}, M)$ by simultaneously
permuting the factors in $\ot_i L_i$ and $H^{\ot I}$.  This is
well-defined because of cocommutativity.

Examples of explicit calculations in $\M^*(H)$ will be provided below.
\end{examples}

A {\em Lie algebra} in a pseudotensor category $\M$ is an object
$L$ together with a polylinear map $\beta\in\Lin(\{L,L\}, L)$ that
satisfies analogues of skew-commutativity and the Jacobi identity:
$\beta=-(12)\beta$, where $(12)\in S_2$, and
$\beta(\beta(\cdot,\cdot),\cdot)=\beta(\cdot,\beta(\cdot,\cdot))-
(12)\beta(\cdot ,\beta(\cdot,\cdot ))$, where $(12)$ now lies in $S_3$.

An {\em associative algebra} is an object $R\in\M$ together with a
polylinear map $\mu\in\Lin(\{R,R\},R)$ satisfying associativity
$\mu(\mu(\cdot,\cdot),\cdot)=\mu(\cdot,\mu(\cdot,\cdot))$.

In $\mathcal{V}ec$ associative (Lie) algebras are just the associative
(Lie) algebras in their usual sense. To avoid confusion we will sometimes
call these algebras {\em ordinary}. The same is true of representations
(to be defined below), cohomology (see \cite{BDK}), etc.

\begin{remark} Consider an associative algebra $(R,\mu)$ in a pseudotensor
category $\M$.  Then the pair $(R, \mu-(12)\mu)$ is a Lie algebra in $\M$ 
(\cite[Prop. 3.11]{BDK}).  However, not every Lie algebra can be
represented in this form, i.e. the PBW theorem does not necessarily hold
in a generic pseudotensor category (see \cite{Ro} for the case of
conformal algebras).
\end{remark}

A {\em representation} of an associative algebra $(R,\mu)$ in $\M$ is an
object $M$ (a {\em module}) together with $\rho\in\Lin(\{R,M\},R)$
satisfying $\rho(\mu(\cdot,\cdot),\cdot)=\rho(\cdot,\rho(\cdot,\cdot))$.
Representations of Lie algebras are defined in the similar way.

\subsection{Preliminaries on Hopf algebras} 
Before proceeding further, we to recall several facts about Hopf
algebras (see \cite[Chapter 1]{Jo} or \cite{Sw} for basic definitions and
notations).

In this paper $H$ will always stand for a Hopf algebra with a coproduct
$\De$, a counit $\eps$, and an antipode $S$. As usual, we use
Sweedler's notations: $\De(h)=h_{(1)}\ot h_{(2)}$ (summation is implied),
$(\De\ot\id)\De(h)=h_{(1)}\ot h_{(2)}\ot h_{(3)}$,
$(S\ot\id)\De(h)=h_{(-1)}\ot h_{(2)}$ etc.

The following formulas will be quite useful:
\begin{align}
\eps(h_{(1)})h_{(2)}&=h_{(1)}\eps(h_{(2)})=h,\lbb{eq.eps}\\
h_{(-1)}h_{(2)}\ot h_{(3)}&=h_{(1)}h_{(-2)}\ot h_{(3)}=1\ot h.\lbb{eq.h}
\end{align}

$H$ also possesses the opposite coproduct $\De^{cop}:h\mapsto h_{(2)}\ot
h_{(1)}$; we denote the Hopf algebra $(H,\De^{cop},\eps,S)$ as $H^{cop}$.
As usual, $H^{op}$ will stand for the algebra $H$ with the opposite
multiplication.

An associative algebra $A$ is called {\em $H$-differential} if it is a
left $H$-module such that
\begin{equation}\lbb{eq.diffalgebra}
h(xy)=(h_{(1)}x)(h_{(2)}y).
\end{equation}

\begin{remark}\lbb{hdiffaspseudo} An $H$-differential algebra is an
associative algebra in the pseudotensor category $\M^l(H)$.
\end{remark}

For an $H$-differential algebra $A$ one defines a {\em smash product}
$A\smash H$ as a tensor product $A\ot H$ of underlying vector spaces with
a new multiplication
\begin{equation*}
(a\smash g)(b\smash h)=a(g_{(1)}b)\smash g_{(2)}h.
\end{equation*}

$G(H)$ stands for the set of group-like elements of $H$,
i.e. $h\in H$ such that $\De(h)=(h\ot h)$.
$P(H)$ is the set of primitive elements, i.e. $h\in H$ with
$\De(h)=1\ot h+h\ot 1$.  Group-like elements form a group with
multiplication inherited from $H$, and $P(H)$ is a Lie subalgebra of $H$
with respect to the standard commutator $[g,h]=gh-hg$.

A typical example of a smash product is $U(P(H))\smash \CC[G(H)]$, where 
$G(H)$ acts on $P(H)$ by inner automorphisms: $g(h)=gpg^{-1}\in P(H)$.

\begin{remark}\lbb{kostant} A theorem due to Kostant \cite[Theorem
8.1.5]{Sw} states that  a cocommutative Hopf algebra $H$ is, in fact,
isomorphic to $U(P(H))\smash\CC[G(H)]$.
\end{remark}

We will also require a standard filtration on $H$:
$\fil^0 H=\CC[G(H)]$ and for $n>0$, 
$\fil^n H =\left\{h\in H|\,\De(h)\in\fil^0 H\ot h+h\ot\fil^0 H+
\sum_{i=1}^{n-1}\fil^i H\ot\fil^{n-i} H\right\}$.

When $H=U(\gg)$ is a universal enveloping algebra, we get the canonical
filtration.  Remark that when $\gg$ is finite-dimensional, $\dim\fil^n
H<\infty$ for all $n$. It is clear that operations on $H$ respect the 
filtration.

When  $H$ is cocommutative, \reref{kostant} implies $\bigcup_n \fil^n
H=H$.  If so, we say that a nonzero element $h\in H$ has {\em degree} n if
$h\in\fil^n H\backslash \fil^{n-1}H$.

In order to define certain operations on pseudoalgebras (see below), we
will need the following:
\begin{lemma}[{\cite[Lemma~2.5]{BDK}}]\lbb{fourier} Every element of $H\ot
H$ can be uniquely represented in the form $\sum_i (h_i\ot 1)\De(l_i)$,
where $\{h_i\}$ is a fixed basis of $H$ and $l_i\in H$.  Also, for any
$H$-module $V$,
\begin{equation*}
(\fil^n H\ot\CC)\De(H)=\fil^n(H\ot H)\De(H)=(\CC\ot\fil^n H)\De(H),
\end{equation*}
where $\fil^n(H\ot H)=\sum_{i+j=n}\fil^i H\ot\fil^j H$.
\end{lemma}

\subsection{Dual algebra of a Hopf algebra}
Denote the dual algebra of $H$ by $X=H^*=\Hom_\CC(H,\CC)$.

It is an $H$-differential algebra with the action defined by 
\begin{equation*}
\langle hx,f\rangle=\langle x,S(h)f\rangle,\text{for $f,h\in H, x\in X$}.
\end{equation*}
One can similarly define the structure of a right $H$-module on
$X$; this makes $X$ into an $H$-bimodule.

$X$ possesses a standard filtration $X=\fil_{-1}X\supset\fil_0 
X\supset\dots$, where $\fil_n X=(\fil^n H)^\perp$. 

When $H$ is cocommutative, $X$ is commutative and $\bigcap_n\fil_n X=0$.

By a basis of $X$ we will always mean a topological basis $\{x_i\}$ such
that for any $n$ only a finite number of $x_i$'s does not lie in $\fil_n
X$ (i.e. $x_i\to 0$ in the standard topology).  Let $\{h_i\}$ be a basis
of $H$ compatible with the standard filtration.  If $\dim\fil^n H<\infty$
for all $n$, the dual basis of $X$ is topological.  For $h\in H$ and $x\in 
X$, we have
\begin{equation*}
h=\sum_i\langle h,x_i\rangle h_i,\quad x=\sum_i\langle x,h_i\rangle x_i,
\end{equation*}
where the first sum is finite and the second converges in the standard
topology.

Define the antipode $S$ on $X$ as a dual of that of $H$: $\langle 
S(x),h\rangle=\langle x,S(h)\rangle$. Also, we introduce a coproduct 
$\De:X\to X\hat\ot X$ where $X\hat\ot X=(H\ot H)^*$ is the completed 
tensor product. By definition, for $x,y\in X$ and $f,g\in H$
\begin{align}
\langle xy,f\rangle&=\langle x\ot y,\De(f)\rangle=\langle x,f_{(1)}\rangle
\langle y,f_{(2)}\rangle,\lbb{eq.deltaf}\\
\langle x,fg\rangle&=\langle \De(x),f\ot g\rangle=\langle
x_{(1)}, f\rangle \langle x_{(2)}, g\rangle.\lbb{eq.deltax}
\end{align}

\begin{remark}\lbb{xactsonh1} For $X$ such that $\dim\fil_n X<\infty$ for
all $n$, one can endow $H$ with the structure of an $X$-differential
algebra.  As in the case of $H$-action on $X$, the action is defined by
$\langle y,xh\rangle=\langle S(x)y,h\rangle$.  The proof is also similar
to that for the $H$-action on $X$ and uses (\ref{eq.deltax}) instead of
(\ref{eq.deltaf}) and (\ref{eq.h}).  Formula (\ref{eq.diffalgebra}) makes
sense as the right-hand side will be finite for every pair of elements of
$H$.
\end{remark}

\subsection{Notations for universal enveloping algebras}
In this paper we mostly restrict our attention to universal enveloping
algebras of finite-dimensional Lie algebras.  Some notations are in order:

Put $H=U(\gg)$ where $\gg$ is a $n$-dimensional Lie algebra spanned over
$\CC$ by $\{\d_i\}^n_{i=1}$.  We fix the canonical (but not PBW) basis of
$H$ indexed by elements of $\ZZ_{\geq 0}^n$:
\begin{equation*}
\d^I=\frac{\d_1^{i_1}\dotsm\d_n^{i_n}}{i_1!\dotsm i_n!},\ \ \ \text{for
}I=(i_1,\dots,i_n).
\end{equation*}
It is easy to see that $\De(\d^I)=\sum_{J+K=I}\d^J\ot\d^K$.

\begin{remark}\lbb{rem.indexsets} For future reference, we need to
describe our notations for the multiindex set $\ZZ_{\geq 0}^n$. The
addition is pointwise. There is the standard partial ordering, i.e.,
$(i_1,\dots,i_n)>(j_1,\dots,j_n)$ iff $i_m>j_m$ for all $m$; if neither
$I\geq J$ nor $J\geq I$, we call $I$ and $J$ {\em incompatible}.  The
index $(0,\dots,0)$ is denoted simply by $0$.  Also, for a multiindex
$I=(i_1,\dots,i_n)$ we put $|I|=i_1+\dots+i_n$ and $(-1)^I=(-1)^{|I|}$.
\end{remark}

The dual Hopf algebra of $H$ is $X=\CC[[t_1,\dots,t_n]]$ with the
canonical dual basis $t^I=t_1^{i_1}\dotsm t_n^{i_n}$. As usual, $t^0=1$.  
The action of $H$ on $X$ is given by differential operators:
$\d_i=-\d/\d t_i$.  The right action is the same: $xh=hx$ for $x\in 
X$, $h\in H$.  Similarly, the action of $X$ on $H$ is defined by 
$t_i=-\d/\d\d_i$.  

\begin{remark}\lbb{filtrwarning}  The standard filtration on $X$ 
regarded as the
dual algebra of $H$ is not the standard decreasing filtration on the
polynomial algebra $\CC[t_1,\dots,t_n]$ (it is shifted by $1$).  In
particular, the standard total degree function on $X$ does not
respect multiplication.
\end{remark}

It is easy to see that $\De(t_i)=1\ot t_i+t_i\ot 1+$summands with both
sides of degree higher than $0$.  This can be generalized for any $t^I$:
\begin{equation}\lbb{coprodform}
\De(t^I)=\sum_{J\leq I} t^J\ot t^{I-J}+\sum_j c_j t^{K_j}\ot t^{L_j},\quad
|K_j|+|L_j|\geq|I|+1,\ c_j\in\CC.
\end{equation}
In particular, one can deduce from (\ref{coprodform}) that
\begin{equation}\lbb{dualfilt}
\De(\fil_{n-1} X)\subset\textsum_{i=0}^n \fil_{i-1} X\hat\ot
\fil_{n-i-1} X.
\end{equation}

\subsection{Pseudoalgebras and their representations}  Recall the
description of the pseudotensor category $\M^*(H)$ (see
\exref{ex.pseudotensor}(3)).

\begin{definition}\lbb{pseudodef} An associative (Lie) {\em pseudoalgebra}
over a cocommutative Hopf algebra $H$ is an associative (Lie) algebra in
$\M^*(H)$.  A representation of a pseudoalgebra a ({\em pseudomodule}) is
its representation in $\M^*(H)$.
\end{definition}

We denote multiplication in a pseudoalgebra $R$ by $*$: $R\ot R\to (H\ot
H)\ot_H R$ and call $a*b$ the {\em pseudoproduct} of $a$ and $b$.

This operation satisfies $H$-bilinearity: for $f,g\in H$, 
$(fa)*(gb)=((f\ot
g)\ot_H 1)(a*b)$, and associativity $(a*b)*c=a*(b*c)$. The explicit
expressions for the latter equality are calculated below, following 
\cite[Chapter 3]{BDK}.

Let $a,b,c\in R$.  To calculate $(a*b)*c\in H^{\ot 3}\ot_H R$ in
accordance with (\ref{pseudotensorhopfcomp}), notice that here
$J=\{1,2,3\}, I=\{1,2\}, \psi_1=\phi=*$, $\psi_2=\id$, and the map $\pi:
J\twoheadrightarrow I$ is given by $\pi(1)=\pi(2)=1, \pi(3)=2$. Put
\begin{align*}
a*b&=\textsum_i (f_i\ot g_i)\ot_H d_i,\\ 
d_i*c&=\textsum_j (f_{ij}\ot
g_{ij})\ot_H d_{ij}.
\end{align*}
Then, as $\De^{(\pi)}=\De\ot\id$, 
\begin{equation}\lbb{pseudoass1}
(a*b)*c=\textsum_{i,j} (f_i{f_{ij}}_{(1)}\ot g_i{f_{ij}}_{(2)}\ot
g_{ij})\ot_H d_{ij}.
\end{equation}
Similarly, for the product $a*(b*c)$, $\De^{(\pi)}=\id\ot\De$, and we
obtain
\begin{align}
a*(b*c)&=\textsum_{ij} (h_{ij}\ot h_i{k_{ij}}_{(1)}\ot
k_j{k_{ij}}_{(2)})\ot_H e_{ij},\lbb{pseudoass2}\\
\text{where}\ &b*c=\textsum_i (h_i\ot k_i)\ot_H e_i,
\ a*e_i=\textsum_j (h_{ij}\ot k_{ij})\ot_h e_{ij}.\notag
\end{align}

For a representation $V$ of $R$, we will also denote the action by $*$:
$a*v\in (H\ot H)\ot_H V$.  It also satisfies $H$-bilinearity and
associativity.

\begin{remark}\lbb{groupaction} Recall that a cocommutative Hopf is a
smash product of $G(H)$ and $U(P(H))$ (\reref{kostant}). For brevity
denote $\Gamma=G(H)$, $H'=U(P(H))$. The action of $\Gamma$ on $H'$ can be
extended to the action on $(H')^{\ot I}$ via $\De^{(I)}(g)=\ot_i g$, and
in an obvious way to the action on $(H')^{\ot I}\ot_{H'} M$ for an
$H'$-module $M$.  It can be shown that the category $\M^*(H)$ is
equivalent to a subcategory of $\M^*(H')$ that consists of $H$-modules and
polylinear maps that commute with the action of $\Gamma$.  Thus, it
follows (\cite[Corollary 5.3]{BDK}) that an $H$-pseudoalgebra is an
$H'$-pseudoalgebra with an action of $\Gamma$ such that $ga*gb=g(a*b)$.
Moreover, the pseudoproduct over $H$ is defined as
\begin{equation}\lbb{eq.groupact}
a*b=\sum_{g\in\Gamma}((g^{-1}\ot 1)\ot_H 1)(ga*b),
\end{equation}
where the products on the RHS are taken over $H'$ and the sum is finite.

This shows that the case of pseudoalgebras over general cocommutative
algebras in most cases reduces to the study of pseudoalgebras over
universal enveloping algebras. 
\end{remark}

\subsection{Annihilation algebra and $x$-products}
As before, $X=H^*$.  Let $R$ be a left module of a cocommutative Hopf
algebra $H$. Define another $H$-module $\coef(R)=X\ot_H R$ with an obvious
left action $h(x\ot_H a)=hx\ot_H a$.  If $R$ is also an
associative $H$-pseudoalgebra, $\coef(R)$ is an associative algebra with
multiplication defined by (\ref{eq.annihmult}).

$\coef(R)$ is an $H$-differential algebra.  Remark also that an 
$R$-pseudomodule $V$ gives rise to an $\coef(R)$-module $\coef(V)=X\ot_H 
V$.

The algebra $\coef(R)$ is called the {\em annihilation algebra} of the
pseudoalgebra $R$.  Its elements $x\ot_H a$ are denoted as $a_x$ and are
called {\em Fourier coefficients} of $a$.

The annihilation algebra closely mirrors the properties of the
corresponding pseudoalgebra.  In particular, when $R$ is
torsion-free, $\coef(R)$ ``distinguishes'' its elements:

\begin{lemma}[cf. {\cite[Proposition 11.5]{BDK}}]\lbb{coeffunique}
Let $M$ be a left $H$-module over a universal enveloping algebra
$H$. All Fourier coefficients of $a\in M$ are zero if and only if $a$ is
torsion.
\end{lemma}

Moreover, sometimes it is possible to get back from $\coef(M)$ to $M$.  
Given a topological left $H$-module $L$, one can construct another module
$\mathcal{C}(L)=\Hom^{\rm cont}_H(X,L)$ (``cont'' stands for continuous in
the standard topologies of $X$ and $H$). Define the map
$\Phi:M\to\mathcal{C}(\coef(M))$ as $\Phi(a)(x)=x\ot_H a$.  In most
interesting cases, $M$ imbeds into $\Phi(M)$ and, if $M$ possesses a
pseudoalgebra structure, so does $\Phi(M)$ (see \leref{CApseudoalg}).

Let $R$ be an associative pseudoalgebra with the pseudoproduct
$a*b=\textsum_i (f_i\ot g_i)\ot_H e_i$.  The choice of $f_i, g_i$, and
$e_i$ is certainly not unique.  By \leref{fourier} we can assume that
$g_i=1$.  This defines the new operation $R\ot R\to H\ot R$: $a\cdot
b=\textsum_i f_i\ot e_i$.  For any $x\in X$ we introduce the {\em
$x$-product}:
\begin{equation}\lbb{xproduct}
a_xb=(\langle S(x),\cdot\rangle\ot\id)a\cdot b=\textsum_i \langle
S(x),f_i\rangle e_i.
\end{equation}
Given the $x$-products of $a$ and $b$, one can also pass back to their
pseudoproduct, obtaining (\ref{eq.pseudomult}).

Notice that the sum in (\ref{eq.pseudomult}) is finite, i.e. for almost
all $x_i$, $a_{x_i}b=0$.

Thus, one can define an associative $H$-pseudoalgebra as a left $H$-module
$R$ equipped with the $x$-products satisfying:
\begin{description}
\item[Locality] 
\begin{equation}\lbb{locality}
\codim\{x\in X\,|\, a_xb=0\}<\infty\quad\text{for any $a,b\in R$};
\end{equation}
\item[$H$-sesquilinearity] 
\begin{equation}\lbb{sesquilin}
\begin{split}
(ha)_xb&=a_{xh}b,\\
a_x(hb)&=h_{(2)}(a_{h_{(-1)}x}b)\quad\text{for any $a,b\in R$, $h\in H$};
\end{split}
\end{equation}
\item[Associativity] 
\begin{equation}\lbb{assoc1}
a_x(b_yc)=(a_{x_{(2)}}b)_{x_{(1)}y}c.
\end{equation}
\end{description}

Locality suggests the following definition: we will call $x\in X$ {\em
maximal with respect to $a$ and $b$} if $a_x b\neq 0$ but for any $y\in
\fil_0 X$ , $a_{xy}b=0$.

Associativity (\ref{assoc1}) can be equivalently stated as
\begin{equation}\lbb{assoc2}
(a_xb)_yc=a_{x_{(2)}}(b_{x_{(-1)}y}c).
\end{equation}

Most of the above properties survive the passing to 
$\mathcal{C}(\coef(R))$:

\begin{lemma}[{cf. \cite[Proposition 11.2]{BDK}}]\lbb{CApseudoalg} Let $R$
be an associative pseudoalgebra. Then $\Phi(R)=\mathcal{C}(\coef(R))$
satisfies (\ref{sesquilin}-\ref{assoc2}).
\end{lemma}

The above formulas and statements, of course, remain true for the
action of $R$ on its pseudomodule $M$ and of $\coef(R)$ on $\coef(M)$.

Using (\ref{eq.pseudomult}) we can obtain formulas similar to
(\ref{assoc1}) and (\ref{assoc2}) for the multiplication in $\coef(R)$:
\begin{equation}\lbb{coefmult}
\begin{split}
a_x\cdot b_y&=(a_{x_{(2)}}b)_{x_{(1)}y},\\
(a_xb)_y&=(a_{x_{(2)}})\cdot(b_{x_{(-1)}y}).
\end{split}
\end{equation}

We can now define structural concepts for associative pseudoalgebras.  
Denote by $A_xB$ the set $\{a_xb\,|\, a\in A, b\in B\}$.  An {\em ideal}
$I$ of $R$ is a pseudoalgebra such that for all $x\in X$, $I_xR\subset I,
R_x I\subset I$.  A pseudoalgebra $R$ whose only ideals are $0$ and $R$ is
called {\em simple}.  A pseudoalgebra $R$ such that for a fixed $n$
$R_{x_1}R_{x_2}\cdots_{x_n}R=0$ for any collection of
$\{x_1,\dots,x_n\}\subset X$ is called {\em nilpotent} (as
$x_i$'s are arbitrary, we can omit the brackets).  A pseudoalgebra that
contains no nilpotent ideals is called {\em semisimple}.  Similar
definitions for pseudomodules will be provided in \seref{sec.rep}.


\section{Examples of Associative Pseudoalgebras}\lbb{sec.ex}

In these section we provide several important examples of associative
pseudoalgebras.  In general, we do not assume that $H$ is a universal
enveloping algebra of a Lie algebra.

\subsection{General examples}
\begin{example}\lbb{ex.cur} Let $H'$ be a Hopf subalgebra of $H$, and $R$
an $H'$-pseudoalgebra.

\begin{definition}\lbb{def.cur} The {\em current extension} of
$R$ is the $H$-pseudoalgebra $\cur_{H'}^H R$ which is the $H$-module 
$H\ot_{H'}R$ with the pseudoproduct $*$ extending the pseudoproduct of $R$ 
by $H$-bilinearity. 
\end{definition}

More explicitly, for $a,b\in R$ with $a*b=\sum_i (f_i\ot
g_i)\ot_{H'}c_i$ we define
\begin{equation}\lbb{eq.cur}
\begin{split}
(f\ot_{H'}a)*(g\ot_{H'}b)&=((f\ot g)\ot_H 1)(a*b)\\
&=\sum_i ((ff_i)\ot(gg_i))\ot_H (1\ot_{H'}c_i)
\end{split}
\end{equation}

The same construction, of course, is possible for any variety of
(pseudo)algebras defined over $H$.

\begin{remark}\lbb{rem.curr}  This terminology is different from that of
\cite{BDK} where current extensions were called current pseudoalgebras.
However, here this term is restricted to a smaller class (see immediately
below).  In the author's view this makes some statements below less
cluttered.
\end{remark}

In particular, when $H'=\CC$, an associative $H'$-pseudoalgebra $R$ is an
associative $\CC$-algebra with the ordinary product.  Then $\cur_{H'}^H
R$, which we will denote simply $\cur R$, has the pseudoproduct
\begin{equation*} 
(f\ot a)*(g\ot b)=(f\ot g)\ot_H (1\ot ab).
\end{equation*}
We will call such a pseudoalgebra $\cur R$ a {\em current pseudoalgebra}.

Current pseudoalgebras have a simple characterization: if in a
pseudoalgebra $R$ $a_xb=0$ for all $a,b\in R$ and $x\in\fil_0 X$, then $R$
is current.

\end{example}

\begin{example}\lbb{ex.rk1}

Let $H=U(\gg)$ be a universal enveloping algebra and $R$ an associative
pseudoalgebra over $H$ that is free and of $\rk 1$ as an $H$-module.
Below we classify such pseudoalgebras.

\begin{lemma}\lbb{rk1} Let $R$ be as above.  Then either the
multiplication in $R$ is trivial (i.e. $a*b=0$ for any $a,b\in R$) or
$R\cong\cur\CC$.
\end{lemma}

\begin{proof} Let $e$ be a generator of $R$ over $H$, i.e. $R=He$.  By
$H$-bilinearity, multiplication in $R$ is completely determined by the
values of the coefficients in the product $e*e$.  Namely, put
\begin{equation*}
e*e=\alpha\ot_H e, \text{\ \ where }
\alpha=\sum_{(I,J)}c_{IJ}\d^I\ot\d^J\in
H\ot H.
\end{equation*}
Then, to classify pseudoalgebras of $\rk 1$, it suffices to classify all
appropriate $\alpha$'s.

Associativity implies
\begin{equation*}
(\alpha\ot 1)(\De\ot\id)(\alpha)=(1\ot\alpha)(\id\ot\De)(\alpha),
\end{equation*}
which can be rewritten as
\begin{equation}\lbb{eq.rk1}
\sum_{(I,J), K+L=I} c_{IJ}\d^I\d^K\ot\d^J\d^L\ot\d^J=
\sum_{(I,J), M+N=J} c_{IJ}\d^I\ot\d^I\d^M\ot\d^J\d^N.
\end{equation}

Pick $I$ with a maximal degree among all such that $c_{IJ}\neq 0$.  Then
by comparing the degrees of the first terms in (\ref{eq.rk1}), we see that
$|K|=0$, thus $|I|=0$.  Similarly, $|J|=0$, i.e. $\alpha=c\ot 1$ for some
$c\in\CC$.  For a non-zero $c$ we can normalize $e$, so that $e*e=(1\ot
1)\ot_H e$.  This makes $R$ isomorphic to $\cur\CC$.
\end{proof}

Lie pseudoalgebras of $\rk 1$ were classified in \cite[4.3]{BDK} by
essentially similar methods.

\end{example}

\subsection{Conformal algebras}
\begin{example}\lbb{ex.conf}  Below we shall partly follow the
introduction to \cite{BDK}.

Let $H=\CC[\d]$ (and $X=\CC[[t]]$).  Then $H$-pseudoalgebras (of any 
variety) are conformal algebras \cite[Chapter 2]{K1}.  A 
non-axiomatic description of such  objects is possible.

Namely, for an arbitrary ordinary algebra, consider the algebra
$A[[z,z^{-1}]]$ whose elements are called formal distribution.  We
introduce bilinear products ${}_{(n)}$, $n\in\ZZ_{\geq 0}$ defined as
\begin{equation*}\lbb{eq.conf}
{f(z)}_{(n)}{g(z)}= Res_{w=0} f(w)g(z)(w-z)^n,
\end{equation*}
where $Res_{w=0}: A[[w,z,w^{-1},z^{-1}]]\to A[[z,z^{-1}]]$ maps a formal
distribution $h(w,z)$ in two variables to the coefficient at $w^{-1}$, and
the product of two formal distributions in different variables is defined
in the usual way.  Clearly, $f(w)g(z)$ can be expressed via the products
${f(z)}_{(n)}{g(z)}$; this is called the operator product expansion.

There is a natural action of $\d=\d_z$ on $A[[z,z^{-1}]]$.  We
call an algebra of formal distributions {\em conformal} if it satisfies
the locality property: for any $f,g$ only a finite number of products $
f_{(n)}g$ is non-zero.  

Every pseudoalgebra $R$ over $H$ can be described in this way for 
$A=X[[t^{-1}]]\ot_H R$.  Operations are related as $f_{(n)}g=f_{t^n}g$.
\end{example}

The above construction obviously generalizes for the case of abelian 
$\gg$.

\subsection{Pseudolinear algebras} Let $V, W$ be $H$-modules. An {\em
$H$-pseudolinear} map from $V$ to $W$ is a linear map $\phi: V\to(H\ot
H)\ot_H W$ such that $\phi(hv)=((1\ot h)\ot 1)\phi(v)$ for $h\in H, v\in
V$.  The space of all such maps, denoted $\chom(V, W)$, is a left
$H$-module: put $(h\phi)(v)=((h\ot 1)\ot_H 1)\phi(v)$. When $V=W$, we
denote the set of all pseudolinear maps as $\cend(V)$.  Though it is
possible to define the action of the product $\phi*\psi$ on $V$ for
$\phi,\psi\in\cend V$, it might not be represented by a finite sum, i.e.,
$\cend V$ is not necessarily a pseudoalgebra.  However, when $V$ is finite
over $H$, $\cend V$ becomes an associative $H$-pseudoalgebra with a
naturally defined multiplication.

\begin{example}\lbb{ex.cend} 
If $V$ is a finite free $H$-module, i.e., $V=H\ot V_0$
for some finite dimensional vector space $V_0$ with a trivial action of
$H$, then $\cend V=H\ot H\ot \End V_0$ with the pseudoproduct defined as
\begin{equation}\lbb{eq.cend}
(f\ot a\ot A)* (g\ot b\ot B)=(f\ot ga_{(1)})\ot_H(1\ot ba_{(2)}\ot AB).
\end{equation}  
(see \cite[Propositions 10.5, 10.11]{BDK}).

Clearly, in the above case $\cend V$ depends only on $\rk V$.  To
emphasize this, for a module of rank $n$ its pseudoalgebra of
endomorphisms will be denoted simply $\cend_n$.

It is not difficult to see that $\cend_n$ is simple
(\cite[Proposition 13.34]{BDK}); however, unlike the case of ordinary
algebras of linear endomorphisms, it is not finite as an $H$-module.

When $H=\CC[\d]$, the conformal algebra $\cend_n$ is sometimes called the
conformal Weyl algebra $\WW_n$ (\cite{Re1}).  In this case the
standard model of a finite $H$-module of $\rk n$ is the module
$E_n=\{a(z)=\sum at^nz^{-n-1}\}_{a\in \CC^n}$.  Here the generators of
$\cend_n$ over $H$ are elements $J^m_A= \sum At^n(-\d_t)^mz^{-n-1}$ where
$A\in\End_n(\CC)$ (\cite[2.?]{K1}).

\end{example}

\subsection{Differential algebras}

\begin{example}\lbb{ex.diff}

Recall that the bialgebra $X^{cop}$ is isomorphic to $X$ as an algebra and 
has the comultiplication $\De^{op}:x\mapsto x_{(2)}\ot x_{(1)}$.  Consider 
an $X^{cop}$-differential algebra $A$, i.e. a topological associative 
algebra with a left $X^{cop}$-action such that for $x\in X^{cop}, a,b\in A$
\begin{equation}\lbb{eq.xcop}
x(ab)=(x_{(2)}a)(x_{(1)}b).
\end{equation}
Recall that $\De(x)\in X\hat\ot X$ is not, in general, a finite sum.
Thus, in order for (\ref{eq.xcop}) to make sense,  we must require that
for any $a\in A$, $\codim\ann a<\infty$. 

\begin{remark}\lbb{xcopconvention} For brevity, the above property will
never be stated in the further exposition but will always be assumed when
we discuss $X^{cop}$-differential algebras. We will simplify terminology
even further and simply call such algebras $X^{cop}$-algebras, always
implying the structure from (\ref{eq.xcop}).
\end{remark}

\begin{remark}\lbb{xactsonh2} For $X$ such that $\dim\fil_nX<\infty$ for
all $n$, a typical example of an $X^{cop}$-algebra is $H^{op}$.  This
statement is ``dual'' to \reref{xactsonh1} and can be deduced in the same
way.
\end{remark}

We introduce the pseudoalgebra structure on
$\dif A=H\ot A$.  Notice that by $H$-sesquilinearity it is enough to
define the products between elements of the type $1\ot a$:
\begin{equation}\lbb{eq.diff}
(1\ot a)_x(1\ot b)=1\ot (a x(b)),\quad\text{for any }x\in X.
\end{equation}
Associativity of these products follows from (\ref{assoc1})
and (\ref{eq.xcop}).  Finite codimension of the annihilator of every $a\in 
A$ implies locality.

Notice that $\dif A$ is generated over $H$ by $1\ot A$.  For brevity we
will denote such elements $1\ot a$ by $\ti a$.

For a free finite $H$-module $V$, $\cend V$ is a differential
pseudoalgebra $\dif H^{op}\ot \End_N(\CC)$.  This may be shown directly
but in the case of arbitrary $H$ the calculations are cumbersome; this
result will follow from a more general statement in the next section (see
\thref{unitisdiff} and \coref{cenddiff}).

However, the case of $H=\CC[\d]$ is much simpler.  Indeed, put
$A=\CC[\d_t]\ot\End_n(\CC)$ and let the generator $t$ of $X=\CC[[t]]$ act
as a derivation in $\d_t$. For $a\in A$, put $\ti a=\sum at^n z^{-n-1}$.
Then $J^m_A=(-1)^m\sum_j (-1)^j \d^j(\wti{a^{(j)}})$ where $a^{(j)}$ is
the $j$-th derivative of $a$ with respect to $\d_t$.  Similarly, $\ti a$
can be expressed via $J^m_A$'s.  This shows that $\cend_n=\dif A$.

\begin{remark}\lbb{currentdiff} A current algebra over a differential
pseudoalgebra is itself a differential pseudoalgebra.  Namely, let $H'$ be
a subalgebra of $H$.  Choose a topological basis of
$X=\CC[[t_1,\dots,t_n]]$ such that $X'=\CC[[t_1,\dots,t_r]]=(H')^*$ for
some $r<n$.  Let $A$ be a differential $(X')^{cop}$-algebra.  One can
consider the induced action of $X^{cop}$ on $A$, namely let $t_i, i>r$, 
act on $A$ trivially.  Then $\cur^H_{H'}\dif_{H'}(A)=\dif_H A$.
\end{remark}

A particular case of the above setting is an arbitrary algebra $A$ with a
trivial $X^{cop}$ action.  Then $\dif_H A=\cur A$.

\end{example}

\begin{example}\lbb{ex.weyl} 

Here we present a differential pseudoalgebra that is neither current, nor
isomorphic to $\cend_n$.

Recall that the $n$-th Weyl algebra $A_n$ is generated by
$\{x_i, y_i\}^n_{i=1}$ such that
\begin{equation*}
\begin{split}
&x_ix_j=x_jx_i,\ \ \ y_iy_j=y_jy_i,\\
&x_iy_j-y_jx_i=\delta_{ij}.
\end{split}
\end{equation*}

Let $H=\CC[\d_1,\dots,\d_{2n}]$.  Then $X=\CC[[t_1,\dots,t_{2n}]]$.
Since $X$ is cocommutative, every $X$-differential algebra gives rise to
an $H$-pseudoalgebra.  

To define the action of $X$ on $A_n$ it is enough to describe the action
of each $t_i$ and check that it conforms to the Leibniz rule (i.e. that
$t_i$ is a derivation of $A_n$).  For $1\leq i\leq n$ put $t_i=\d/\d x_i$
and for $n+1\leq i\leq 2n$, $t_i=\d/\d y_i$.  Less formally, we put for
$i\leq n$, $t_i=-\ad y_i$, and for $i>n$, $t_i=\ad x_i$; this immediately
implies the Leibniz rule.

Remark that as $A_n$ is simple, $\dif_H A_n$ is also simple (see
\leref{xcopsimple}).

This examples generalizes to the case of $H=U(\gg)$.  Let $\hat\gg$ be a
one-dimensional central extension of $\gg$: $0\to \CC {\tt
c}\to\hat\gg\to\gg\to 0$.  Let $\phi$ be the corresponding cocycle. The
action of $X$ (or $X^{cop}$) on $\gg$ extends trivially to $\hat\gg$.  
Put $A=U(\hat\gg)/(1-{\tt c})$.  The $X^{cop}$-action on $\hat\gg$ passes 
to $A$.  Then $\dif A$ is a simple $H$-pseudoalgebra (this can be 
demonstrated directly or via \leref{xcopsimple}).  We will denote it 
$\cend_n^\phi$.

\begin{remark}\lbb{trivcocycle} When $\phi$ is the trivial cocycle,
$\cend_n^\phi=\cend_n$.
\end{remark}

For an abelian $\gg$, $\hat\gg$ is the Heisenberg algebra, and $A$ as
described above is the Weyl algebra.
\end{example}


\section{Unital Pseudoalgebras}\lbb{sec.unital}

In this section we will only consider the case of  $H=U(\gg)$ where $\gg$
is a finite-dimensional Lie algebra.  We are interested in some sort of
classification of associative $H$-pseudoalgebras and their
representations.

\subsection{Definition of unital algebras}
Any classification of ordinary algebras begins with that of {\em unital}
algebras.  The trivial observation is that an ordinary algebra $A$ is
unital (i.e. possesses an identity) if there is an embedding $\CC\to A$
that agrees with the $\CC$-action on $A$.

We shall introduce a similar concept for $H$-pseudoalgebras.  The role of
$\CC$ will be played by the ``smallest'' pseudoalgebra $\cur\CC$ (cf. 
\leref{rk1}).  

In order to define unital pseudoalgebras, we shall study in greater detail
the representations of $\cur\CC$.  From now on, we will always denote the
generator of $\cur\CC$ as an $H$-module as $e$.

\begin{lemma}\lbb{curcrep} {\rm (i)} Let $V$ be a $\cur\CC$-module.
Then $V=V^0\oplus V^1$, where $V^0$ and $V^1$ are submodules of
$V$ such that $e*V^0=0$ and for every $v\in V^1, e_1v=v$.

{\rm (ii)} $V^1$ is a torsion-free $H$-module.
\end{lemma}

\begin{proof} (i) $\cur\CC$ is an ordinary associative algebra with
respect to the product ${}_1$, hence $V$ splits into the direct sum of
ordinary submodules $V^0\oplus V^1$ such that $e$ acts as a multiplication
by $i$ on $V^i$.

For any $x\in X$, if $e_1v=0$, we have
$0=e_x(e_1v)=(e_{x_{(2)}}e)_{x_{(1)}}v=(e_1e)_xv=e_xv$, thus $e*V^0=0$.

A direct calculation shows that $V^0$ and $V^1$ are $H$-stable.

(ii)  Assume that $v\in V^1$ is torsion, i.e., there exists $h\in H$ such
that $hv=0$.  Suppose we can choose $x$ such that $S(h)x$ is maximal with
respect to $e$ and $v$.  But since
\begin{equation*}
0=e_x(hv)=h_{(2)}(e_{h_{(-1)}x}v)=e_{S(h)x}v\neq 0,
\end{equation*}
this is impossible and for all $x\in X$, $e_xv=0$.  Hence,$v\in V^0\cap
V^1$ and $v=0$.
\end{proof}

\begin{remark}\lbb{alsofree} In fact, we will show in the
proof of \leref{freerep} that $V^1$ is free as an $H$-module, but for now
torsion-freeness will suffice.
\end{remark}

\begin{definition}\lbb{def.unit} An $H$-pseudoalgebra $R$ is called
{\em unital} if 

1) there exists an embedding $\cur\CC\to R$;

2) as a $\cur\CC$-module $R$ has no zero component $R^0$.
\end{definition}

We shall denote the image of the generator of $\cur\CC$ in $R$ by $e$ as
well and call it the {\em pseudoidentity} of $R$.

Differential pseudoalgebras (see \exref{ex.diff}) over unital algebras are
unital and, since identity can be adjoined to any ordinary algebra, any
differential pseudoalgebra can be embedded into a unital one. Thus,
speaking of differential pseudoalgebras, we will always assume them to be
unital.

\begin{remark}\lbb{embed} It is unknown, in general, what pseudoalgebras
can be embedded into unital ones.  Torsion-freeness over $H$ is a
necessary condition (\leref{unit.free}), and one can provide a number of
sufficient conditions as well, e.g. having a faithful finite
representation (then there is an embedding into $\cend_n$).
\end{remark}

\begin{remark}\lbb{idnetitynotunique} Unlike the case of ordinary algebras
with an ordinary identity, pseudoidentity is not unique. Consider, for
instance, the conformal algebra (i.e. a pseudoalgebra over $\CC[\d]$)
$\cur\End_n(\CC)$. Clearly, $\ti 1$ is a pseudoidentity, but so is $\ti
1+\d\ti r$ for any nilpotent $r$ of nilpotency degree $2$ (i.e. $r^2=0$).
\end{remark}

Nonetheless, unital algebras possess a number of good properties.

\begin{lemma}\lbb{unit.free} Let $R$ be a unital $H$-pseudoalgebra.  Then
$R$ is a torsion-free $H$-module.
\end{lemma}

\begin{proof} \leref{curcrep}, (ii).\end{proof}

\subsection{Classification}
Remark that $e_1$ acts in $\coef(R)$ as a left identity.  It is easy to
construct an example of a pseudoalgebra such that $\coef(R)$  possesses no
right identities. E.g. consider $R=\cur A$ where $A$ has no right
identities, then $\coef(R)$ is a tensor product of algebras $X\ot A$.  In
order to provide a good classification of unital algebras, we need to
exclude such degenerate examples.

\begin{definition}\lbb{defann} The {\em left annihilator} of a 
pseudoalgebra $R$ is the set of elements $a$ such that $a*b=0$ for all 
$b\in R$.
\end{definition}

It is clear that $L(R)$ is an ideal of $R$.

\begin{lemma}\lbb{unitann} For a unital $R$, $L(R)=\{a|\,a*e=0\}$.
\end{lemma}

\begin{proof}  Let $a$ be such that $a*e=0$.  For any $x\in X, b\in R$,
$a_xb=a_x(e_1b)=(a_{x_{(2)}}e)_{x_{(1)}}b=0$.  Hence, $a*b=0$.
\end{proof}

\begin{theorem}\lbb{unitisdiff} A unital pseudoalgebra $R$ with a zero
left annihilator is differential: $R=\dif A$ for some associative $A$.  
Moreover, if $R$ is finitely generated as a pseudoalgebra, $A$ is a 
finitely generated algebra.
\end{theorem}

\begin{proof} Consider the subset $A=1\ot_H R$ of the annihilator algebra
$\coef(R)$.  Clearly, it is a subalgebra of $\coef(R)$ with a left
identity $1\ot e$. We will show that for a unital pseudoalgebra, $R=\dif
A$.

We shall describe, at first, the annihilator subalgebra of $\cur\CC$.
Recall that $e*e=(1\ot 1)\ot_H e$. Since $(x\ot e)(y\ot e)=xy\ot e$, $e_1$
is the identity in $\coef(\cur\CC)$ and $\coef(\cur\CC)$ is generated by
$\{e_{t_i}\}$.

Assume for now that $e_1$ is the left and right identity in $\coef(R)$.
Since $e_x=e_1e_x$, $e_x$ is not a zero divisor.

Suppose $R$ is generated over $H$ by the set $R_0$ of elements $a$ such
that 
\begin{equation}\lbb{eq.goodbasis}
a*e=(1\ot 1)\ot_H a 
\end{equation}
(i.e. $a_1e=a$ and $1\in X$ is maximal with respect to $a$ and $e$).  
Notice that if $a$ is such a generator, then so is $e_xa$ for every
$x\in X$.  

Clearly $a_x=(a_1e)_x=a_1e_x$. As the collection $\{a_x\}$ is unique for
each element of $R$ (see \leref{coeffunique}), we conclude that the above
set of generators of $R$ is in 1-1 correspondence with $A$.  Moreover,
$a_1b_1=(a_1b)_1$, hence $A$ is an algebra with multiplication determined
by the ${}_1$-product in $R$.

Define the action of $X$ on $A$: $x(a_1)=(e_xa)_1$.  The following
calculation implies that $x(a_1b_1)=x_{(2)}(a_1)x_{(1)}(b_1)$:
\begin{equation*}
\begin{split}
e_x(a_1b)&=(e_{x_{(2)}}a)_{x_{(1)}}b=(e_{x_{(2)}}(a_1e))_{x_{(1)}}b\\
&=((e_{x_{(2)}}a)_1e)_{x_{(1)}}b=(e_{x_{(2)}}a)_1(e_{x_{(1)}}b).
\end{split}
\end{equation*}

Therefore, $A$ is an $X^{cop}$-differential algebra.  Hence, as $R$ is
torsion free over $H$, it follows from \leref{coeffunique} that the
$x$-products of elements from $R_0$ satisfy (\ref{eq.diff}).

To show that with the $X^{cop}$ action defined as above, $R\simeq\dif A$,
it remains to prove that $R=H\ot R_0\simeq H\ot A$ as $H$-modules, i.e.
that $R$ is a free $H$-module generated by $R_0$.  Assume the contrary,
namely, that there exist non-zero elements $b_I\in R_0$ such that $\sum_I
\d^I b_I=0$ for some finite collection of $I$'s.  Among these $I$'s,
choose a maximal $J$ (with respect to the natural ordering of $n$-tuples).  
Then $0=(\sum_I \d^I b_I)_{t^J}e=(-1)^J b_J$, a contradiction. Hence,
$R=H\ot R_0$, where $R_0$ is a generating set of $R$ satisfying
(\ref{eq.goodbasis}).  By construction of $A$, $R=\dif A$.

Now it remains to construct such a generating set $R_0$.

Fix an arbitrary element $a\in R$.  For $I$ such that $t^I$ is maximal
with respect to $a$ and $e$, put $a_I=(-1)^I a_{t^I}e$.  A direct
calculation utilizing (\ref{coprodform})  shows that for any $J$,
$(a_I)_{t^J}e=\delta_{0,J}a_I$, i.e. $a_I$ satisfies (\ref{eq.goodbasis}).
Consider now the element $a-\d^Ia_I$.  For $J$ such that $J\not\geq I$ and
$a_{t^J}e=0$, $\d^It^J=0$, thus clearly, $(a-\d^Ia_I)_{t^J}e=0$.  Also,
for $J>I$, $(a-\d^Ia_I)_{t^J}e=-(a_I)_{t^{J-I}}e=0$.  Finally,
$(a-\d^Ia_I)_{t^I}e=a_I-(a_I)_1e=0$ as well.

We conclude that by subtracting from $a$ elements of the type $hb$ where
$h\in H$ and $b$ satisfies (\ref{eq.goodbasis}), we can lower the number
of $I$'s such that $a_{t^I}e\neq 0$.  Since in the process we also lower
the degree of such $t^I$'s, we will at some point obtain an element $c$
such that $c_{t^I}e=0$ for $I>0$.  Then either $c_1e=0$ as well or
$c_1e\neq 0$. In the former case, $c*e=0$, hence $c\in L(R)$ and $c=0$.
In the second case, for arbitrary $d$ and $x$, $c_xd=c_x(e_1d)=(c_1e)_xd$,
hence $c-c_1e\in L(R)$ and we see that $c$ satisfies (\ref{eq.goodbasis})
as well.

Therefore, given a set of $H$-generators of $R$, we can construct a set of
generators satisfying (\ref{eq.goodbasis}) which shows that $R$ is a
differential pseudoalgebra.

Moreover, given a set that generates $R$ as a pseudoalgebra, the elements
obtained from it by the above procedure will also generate $R$ and their
tensor products with $1$ will generate $A$.  Due to locality, if $R$ is
finally generated, we will produce a finite number of generators for $A$.

We turn now to the general case. Let $\mathcal{B}(R)=\coef(R)e_1$.  
Clearly $\mathcal{B}(R)$ is a topological associative algebra.  Since for
any $1\neq h\in H$, $h(1\ot e)=0$, for any $a\in R$ and $x\in X$,
$h(a_xe_1)=h(a_x)1(e_1)=h(a_x)e_1$.  We see that $\mathcal{B}(R)$ is an
$H$-differential algebra as well. Hence, $\bar
R=\mathcal{C}(\mathcal{B}(R))=\Hom^{\rm cont}_H (X,\mathcal{B}(R))$
satisfies (\ref{sesquilin}-\ref{assoc2}) (cf. \leref{CApseudoalg}).

Define a map $\phi:R\to\bar R$, $\phi(a)=a'$, where $a'(x)=a_xe_1$.  By
definition of multiplications in $\bar R$, for $a,b\in R$ and $x,y\in X$,
\begin{equation}\lbb{eq.prim}
\begin{split}
(a'_xb')(y)=&a'(x_{(2)})b'(x_{(-1)}y)=a_{x_{(2)}}e_1b_{x_{(-1)}y}e_1=\\  
=&(a'_xb')(y)=a_{x_{(2)}}b_{x_{(-1)}y}e_1=(a_xb)_ye_1, 
\end{split}
\end{equation}
as $e_1$ is the left identity in $\coef(R)$.  Thus we obtain that
$\phi(a_xb)=a'_xb'$.  Denote $\Im\phi$ by $R'$.  The calculation in
(\ref{eq.prim}) shows that multiplications in $R'$ are also local, hence
$\phi$ is a pseudoalgebra map.  We also conclude that $R'$ is a unital
pseudoalgebra with pseudoidentity $e'$.

The above construction of a generating set satisfying (\ref{eq.goodbasis})
could be repeated for $R'$ with $\coef(R')$ replaced with $\mathcal{B}(R)$
with the conclusion $R'=\dif \mathcal{B}(R)$. Therefore, if $R\simeq R'$,
the proof will be finished.

Assume $\phi$ is not injective.  Thus, there exists $a\in R$ such that
$a_xe_1=0$ for all $x\in X$.  As $e_x=e_1e_x$, $a_xe_y=0$ for all $y\in X$
too.  This implies that $(a_xe)_y=0$ for all $y$, hence $a_xe=0$.
Therefore, by \leref{unitann}, $a=0$.
\end{proof}

The above theorem is especially useful in the cases described below:

\begin{corollary}\lbb{cenddiff} $\cend_n$ is a differential algebra over
$\End_n(\CC)\ot H^{op}$.
\end{corollary}

\begin{proof} Definition (\ref{eq.cend}) of the pseudoproduct in $\cend_n$
implies that it is unital with a zero left annihilator (the latter also
follows from the simplicity of $\cend_n$, see \coref{semisimdiff}
immediately below), hence it is differential. Repeating the construction
of $A$ in the proof of \thref{unitisdiff} in this particular case gives
the description of the underlying $X^{cop}$-algebra.  Namely, we have
\begin{equation*}
A=\left\{1\ot_H(1\ot b\ot B)\,|\,b\in H, B\in\End_n(\CC)\right\} 
\end{equation*}
with the multiplication
\begin{equation*}
\left( 1\ot_H(1\ot b\ot B)\right) \left(1\ot_H(1\ot c\ot C)\right)=
1\ot_H(1\ot cb\ot BC),
\end{equation*}
which shows that $A$ is isomorphic to $\End_n(\CC)\ot H^{op}$.
\end{proof}

\begin{corollary}\lbb{semisimdiff} A unital semisimple pseudoalgebra is
differential.
\end{corollary}

\begin{proof} By associativity the left annihilator is a nilpotent ideal.
\end{proof}

\subsection{Unital pseudoalgebras over cocommutative Hopf algebras} We
briefly discuss here what happens in the case of a more general
cocommutative $H$.
 
Let $R$ be a semisimple unital pseudoalgebra over a cocommutative Hopf
algebra $H$. Recall (\coref{groupaction}) that here $R$ is a pseudoalgebra
over $H'=U(P(H))$.  Interpretation (\ref{eq.groupact}) of pseudoproduct
over $H$ in terms of that over $H'$ together with (\ref{xproduct}) imply
that $R$ is unital over $H'$.  In particular, if $e$ is a pseudoidentity
over $H$, it remains such over $H'$ (the calculation is direct and is,
therefore, omitted).

Consider now the pseudoproduct $(ge)*e=(g\ot 1)\ot_H e$ over $H$.  It does
not survive passing to $H'$ (cf. the construction in \cite[Chapter
5]{BDK}), hence $L(R)\neq 0$ by \leref{unitann}.  Remark that
according to \deref{defann}, it is impossible to find another
pseudoidentity in $R$ regarded as an $H'$-pseudoalgebra such that $R$
would satisfy the conditions of \thref{unitisdiff}.  Hence, $R$ is not an
$H'$-differential pseudoalgebra over a unital algebra.
 
On a more general note, unitality as defined above does not seem to be the
right concept for the study of pseudoalgebras over a generic cocommutative
algebra: for instance, there is no analogue of \leref{rk1}.


\section{Representations of Unital Pseudoalgebras}\lbb{sec.rep}

We turn now to the description of representations of unital
pseudoalgebras.  Although we will work only with unital differential
pseudoalgebras, in light of \thref{unitisdiff}, this simply means that we
impose a technical condition $L(R)=0$.  Since most of the interesting
pseudoalgebras are semisimple, this holds automatically (
\coref{semisimdiff}).

The goal is to provide a statement similar to \thref{unitisdiff}, i.e. to
establish a correspondence between the categories of modules of $A$ and
$\dif A$.  As in the previous section, $H=U(\gg)$ for a finite-dimensional
$\gg$.

\subsection{Structure of representations of unital algebras}
Consider a representation $V$ of a differential pseudoalgebra $R=\dif A$.
Recall that by \leref{curcrep}, as a $\cur\CC$-module,
$V=V^0\oplus V^1$, where $e*V^0=0$.  Thus for any $\ti a, a\in A$ and
$v\in V^0$, $\ti{a}_xv=(\ti{a}_1e)_xv=0$, and $R*V^0=0$.  Also, since
$e_1(a_xv)=(e_1a)_xv=a_xv$, $V^1$ is $R$-stable.  Therefore, the
decomposition of $V$ is valid over $R$ as well.

\begin{definition} A module $V$ of a unital differential
pseudoalgebra $R$ is {\em unitary} if it has no zero component $V^0$.
\end{definition}

Let now $R=\dif A$ be a unital differential pseudoalgebra and $V$ its
unitary module.

As in the proof of \thref{unitisdiff}, for any $v\in V$, we can consider
elements $v_I=e_{t^I}v$.  If $t^I$ is maximal with respect to
$e$ and $v$, $e_{t^J}v_I=\delta_{0,J}v_I$ for any $n$-tuple $J$.  Now,
consider the difference $w=v-\d^Iv_I$.  Direct calculations show that
$e_{t^J}w=0$ for $J$ such that either $J\geq I$ or $J$ is incompatible
with $I$ and $e_{t^J}v=0$.  By taking such differences repeatedly we will
arrive at $w$ such that $e_{t^I}w=\delta_{0,I}w$.  Hence, $V$ is
generated over $H$ by elements $v$ such that $e*v=(1\ot 1)\ot_H v$.

For such an element $v$, $\ti a_{t^I}v=(\ti a_1e)_{t^I}v=\ti
a_1(e_{t^I}v)$, hence 
\begin{equation}\lbb{eq.unitrep}
\ti a*v=(1\ot 1)\ot_H (\ti a_1v).  
\end{equation}

\begin{lemma}\lbb{freerep}
A unitary module of a unital differential pseudoalgebra is free as
an $H$-module.
\end{lemma}

\begin{proof} 
Assume the contrary, i.e. the existence of a finite collection of
non-zero $v_I\in V_0$ such that $\sum_I \d^Iv_I=0$.  Pick $J$ to be a
maximal $n$-tuple such that $v_J\neq 0$.  Then, by calculating
$e_{t^J}(\sum_I \d^Iv_I)$ via (\ref{eq.unitrep}), we see that $v_J=0$.
\end{proof}

\begin{corollary}\lbb{freeunitdiff} A unital differential pseudoalgebra is
free as an $H$-module.
\end{corollary}

Now, since $V$ is free over $H$, we see that $V=H\ot V_0$, and $\ti a$
acts on $V_0$ in accordance with (\ref{eq.unitrep}).

Remark that $\coef(V)=X\ot_H V_0=X\ot V_0$.  In particular, if we
write out $\ti a_1v=v_0+\sum_i h_iv_i$, where $v_i\in V_0$, then $1\ot_H
(\ti a_1v)=1\ot v_0$.  Recall that $A=1\ot_H R$ and $a=1\ot_H\ti a$. We
can introduce the action of $A$ on $V_0$ viewed as the subspace $1\ot
V_0$ of $\coef(V)$.  Thus, 
\begin{equation}\lbb{eq.ordinrep}
av=a(1\ot_Hv)=1\ot_H(\ti a_1v),\ \ \ v\in V_0
\end{equation}

We sum up the above discussion in the following lemma:

\begin{lemma}\lbb{repfacts}  Let $V$ be a module of a unital
differential pseudoalgebra $R=\dif A$.  Then $V=V^0\oplus V^1$ where
$R*V^0=0$ and $V^1=H\ot V_0$ with the action of $R$ on elements of $V_0$
described by (\ref{eq.unitrep}).  Moreover, there is a structure of an
$A$-module on $V_0$ described by (\ref{eq.ordinrep}).
\end{lemma}

\subsection{Constructing representations of unital pseudoalgebras}
Conversely, let $M$ be a unitary left module for an $X^{cop}$-differential
algebra $A$.  Our goal is to construct a related representation of $R=\dif
A$ on a left $H$-module $\ti M=H\ot M$.  Before this, we shall endow $X\ot
M$ with the structure of an $\coef(R)$-module.  (Even though
$\coef(R)=X\ot A$, we will write its basis elements as $x\ot_H\ti a$ to
emphasize the relation with $R$).

Naturally, we put $(1\ot_H\ti a)(1\ot m)=(1\ot am)$ and $(x\ot_H e)(1\ot
m)=(x\ot m)$.  In general, $(x\ot_H\ti a)(y\ot m)=(1\ot_H\ti a)(xy\ot_H
e)(1\ot m)=(xy\ot_H \ti a)(1\ot m)$.  Hence, to describe explicitly the
action of of $\coef(R)$ on $X\ot M$, it suffices to write out the
expression for $(x\ot_H\ti a)(1\ot m)$.  Recall that $x\ot_H\ti a$ is $\ti
a_x=\ti a_1e_x=(e_1\ti a)_x$.  Using (\ref{assoc1}) and (\ref{assoc2}), it
is not difficult to see that 
\begin{equation*}
\begin{split}
\ti a_1e_x=&\sum_I (e_{(\d^Ix_{(2)})t^I}\ti a)_{x_{(1)}}=\\
=&\sum_I (\d^Ie)_x(e_{t^I}\ti a)_1,
\end{split}
\end{equation*}
where the first equality is valid because $\sum_I (\d^It^J)t^I=0$ whenever
$J\neq 0$.

Therefore, $x\ot_H\ti a=\sum_I(\d^Ix\ot_H e)(1\ot_H\wti{t^Ia})$,
and we obtain
\begin{equation}\lbb{eq.repcoeff}
(x\ot_H\ti a)(1\ot m)=\sum_I (\d^Ix)\ot (t^I(a)m).
\end{equation}

We now turn to the description of $\ti M$.  First of all, remark that
$\coef(\ti M)=X\ot M$.  Thus, for $x,y\in X$,
$(e_x\ti m)_y=e_{x_{(2)}}\ti m_{x_{(-1)}y}=(x_{(2)}x_{(-1)}y)\ot
m=\eps(x)y\ot m$ (\ref{eq.eps}).  Hence, if $\eps(x)=0$, $e_xm=0$, and
$e*\ti m=(1\ot 1)\ot_H \ti m$.  That is, $\ti M$ is a unitary module and,
according to (\ref{eq.unitrep}), $\ti a*\ti m=(1\ot 1)\ot_H \ti a_1\ti m$.
It remains only to determine $\ti a_1\ti m$ for arbitrary $a\in A, m\in
M$.  Checking the coefficients, we obtain from (\ref{eq.repcoeff}):
\begin{equation}\lbb{eq.tilderep1}
\ti a_1\ti m=\sum_I \d^I\big(\wti{t^I(a)m}\big).
\end{equation}

We summarize the above discussion as

\begin{lemma}\lbb{tilderep} Let $A$ be a unital $X^{cop}$ differential
algebra, and $M$ an $A$-module.  Then $\ti M=H\ot M$ is a
representation of $\dif A$ with the action described by
\begin{equation}\lbb{eq.tilderep}
\ti a*\ti m=(1\ot 1)\ot_H\Big(\,
\sum_I \d^I\big(\wti{t^I(a)m}\big)\, \Big).
\end{equation}
\end{lemma}

\subsection{Classification and corollaries}
Summing up, we obtain the full description of representations of unital
differential pseudoalgebras.

\begin{theorem}\lbb{unitrep} Let $V$ be a module of a unital
differential pseudoalgebra $R=\dif A$. Then $V=V^0\oplus V^1$, where
$R*V^0=0$ and $V^1=\ti M$ for some $A$-module $M$.  In particular, $V^1$
is free over $H$.
\end{theorem}

\begin{proof} The decomposition of $V$ as well as the $0$ action of 
$R$ on $V^0$ follow from \leref{repfacts}.  Freeness of $V^1$ is explained
in \leref{freerep}, in particular, we know (again, from \leref{repfacts})
that $V^1=H\ot V^0$, where $V^0$ is an $A$-module.  

We can construct another representation $\wti{V_0}$ of $R$.  By comparing
(\ref{eq.unitrep}) with (\ref{eq.tilderep}), we see that the $R$-action on
both $V^1$ and $\wti{V_0}$ is determined by the ${}_1$-action only.  Now,
define the degree of $a\in A$ as the maximal value of $|I|$ such that
$t^I(a)\neq 0$.  Inducting on the degree and comparing (\ref{eq.ordinrep})
with (\ref{eq.tilderep1}), we conclude that $\wti{V_0}\simeq V^1$.
\end{proof}

\begin{remark}\lbb{algisnotrep} The proofs of both \thref{unitisdiff} and
\thref{unitrep} required constructions of particular $H$-generating sets
of, respectively, a given pseudoalgebra and a given module.  However, if
one considers a unital algebra as a $\cur\CC$-module, these bases
are clearly different (e.g., compare (\ref{eq.goodbasis}) and
(\ref{eq.unitrep})). For conformal algebras $\cend_n$ both were written
out explicitly in \exref{ex.diff}.

For the general case of pseudoalgebras $\cend_n$, the basis from
\exref{ex.cend} is the one corresponding to its structure as a
$\cur\CC$-module: $\cend_n=\wti{M_n}$ where $M_n=H\ot\End_n(\CC)$.
\end{remark}

\begin{remark}\lbb{zeromod} Clearly, a non-unitary $A$-module $M$ gives
rise to a non-unitary $\dif A$-module $\ti M=H\ot M$.  However, in this
case the converse is not true.  For example, a non-unitary $\dif
A$-module does not have to be free over $H$.
 
Nonetheless, for consistency we will sometimes use the notation $\ti M$
for non-unitary modules.  In particular, we will denote the
zero-dimensional $\dif A$-module as $\ti 0$.  
\end{remark}

We now turn to the structural theory of representations of unital
algebras; obviously, \thref{unitrep} will be our main tool.

The definitions are the same as in the ordinary case.  A module $V$ over
a pseudoalgebra is called {\em irreducible} if it contains no submodules
except for $0$ and $V$, {\em indecomposable} if it can not be presented as
a sum of two non-zero submodules, and {\em completely reducible} if it
decomposes into a direct sum of irreducible ideals.

\begin{corollary}\lbb{unitsubmod} Let $\ti M$ be a unitary module of a
unital differential pseudoalgebra $R=\dif A$ and $W$ its submodule.  Then
$W=\ti N$ for an $A$-module $N\subset M$.
\end{corollary}

\begin{proof} The argument is the same as in the proof of \leref{freerep}.

Put $N=\{m|\,\ti m\in W\}$.  Obviously, $\ti N\subset W$. \thref{unitrep}
implies that $W$ is unital as well, so we may apply (\ref{eq.unitrep}).  
For an arbitrary element $w=\sum_I \d^I\wti{m_I}\in W, m_I\in M$, let $J$
be a maximal $n$-tuple among $I$'s such that $m_I\neq 0$.  Then
$e_{t^J}w=(-1)^J m_J\in W$. By induction, we obtain that all $m_I$ lie
in $W$, and $w\in\ti N$.
\end{proof}

\begin{corollary}\lbb{irredrep}  Let $V$ be a module over a unital
differential pseudoalgebra $R=\dif A$.  Then $V$ is irreducible if and
only if $V=\ti M$ for an irreducible $A$-module $M$ (not necessarily
non-zero).
\end{corollary}

\begin{proof}  Since $V=V^0\oplus V^1$, either of the components must be
$0$.  If $V=V^1$, then by \thref{unitrep} $V=\ti M$ and, clearly, $M$ must
be irreducible as well.  If $V=V^0$, then every element of $V$ gives rise
to an $R$-submodule, i.e. $V=\ti 0$.

Conversely, by \coref{unitsubmod}, a non-zero irreducible $A$-module $M$
gives rise to an irreducible $R$-module $\ti M$.
\end{proof}

Similarly, we can prove:

\begin{corollary}\lbb{indecrep} Let $V$ be a module over a unital
differential pseudoalgebra $R$.  Then $V$ is indecomposable if and only if
$V=\ti M$ for an indecomposable $A$-module $M$ (not necessarily non-zero).
\end{corollary}

\begin{corollary}\lbb{complred} Let $V$ be a module over a unital
differential pseudoalgebra $R$.  Then $V$ is completely reducible if and
only if $V=\ti M$ for a completely reducible $A$-module $M$.
\end{corollary}

\begin{proof}  Clearly, $V^0$ is completely reducible if and only if
$V^0=\ti 0$.  Complete reducibility of $V^1$ again follows from
\thref{unitrep} and \coref{unitsubmod}.
\end{proof}

\begin{remark}\lbb{faithfulmod} The only major notion that does not
immediately carry over from $A$-modules to $\dif A$-modules if 
faithfulness. Indeed, if $\ti M$ is faithful, $M$ need not be.  The right 
concept here is to require that the annihilator of $M$ does not 
contain any $X^{cop}$-stable ideals.
\end{remark}

\subsection{Representations of pseudolinear conformal algebras}
The above corollaries allow us, for example, to classify the 
representations of the pseudoalgebra $\cend_n$.  Below we shall do so in 
the case of conformal algebras (see \exref{ex.conf}).

We will completely describe irreducible and indecomposable modules
over the conformal algebra $\WW_n=\cend_n$.  By the above corollaries this
comes down to explaining how $\End_n(\CC)\ot\CC[\d]$-modules look like.

In \exref{ex.cend} we described the standard module $E_n$ over $\WW_n$.
By \coref{cenddiff} $\cend_n=\dif A$ where $A$ is the algebra of 
$n\times n$ matrices over $\CC[\d]$. Thus $E_n=\ti M_n$ where the 
$A$-module $M_n$ is an $n$-dimensional vector space on which $\d$ acts as 
the identity operator.

This can be generalized to the case of the module $M_n^\alpha$ which is
again an $n$-dimensional space on which $\d$ acts as
$\alpha\in\End_n(\CC)$.  Thus we obtain a family of modules
$E_n^\alpha=\wti{M_n^\alpha}$ which can be explicitly written as
\begin{equation*}
E_n^\alpha=\left\{a(z)=\sum (at^ne^{-\alpha t})z^{-n-1},\, 
a\in\CC^n\right\},\quad\alpha\in\End_n(\CC).
\end{equation*}

Every irreducible $\CC[\d]$-module is one-dimensional ($M_1^\alpha$ in
the above notations) and every irreducible $\End_n(\CC)\ot\CC[\d]$-module
is of the form $M_n^\alpha$.  Thus, \coref{irredrep} implies

\begin{proposition}\lbb{weylirred} Finite irreducible modules over $\WW_n$
are of the form $E_n^\alpha$ with the standard action.
\end{proposition}

The case of indecomposable modules is similar.  Let $U$ be a
finite-dimensional space with an indecomposable endomorphism
$\alpha$.  Then $U$ is an indecomposable $\CC[\d]$-module
($x\mapsto\alpha$) and every indecomposable $\CC[\d]$-module has this
form.  For $\End_n(\CC)\ot\CC[\d]$ the module $M_n^\alpha(U)=\CC^n\ot U$
with the obviously defined action is indecomposable.  Remark that the
irreducible modules $M_n^\alpha$ defined above are simply
$M_n^\alpha(\CC)$.  

These describe all indecomposable modules over $\End_n(\CC)\ot\CC[\d]$
(consider the decomposition of such as a $\CC[\d]$-module; all components
will be of the same height).  We thus obtain indecomposable $\WW_n$
modules $E_n^\alpha(U)=\wti{M_n^\alpha(U)}$. In \cite{K2} they
were denoted $\sigma_\alpha^{\mathit{as}}$ (``\textit{as}'' stands
for associative); the following result was first stated there as well:

\begin{proposition}\lbb{weylindec} Finite indecomposable modules over
$\WW_n$ are exactly $E_n^\alpha(U)$.
\end{proposition}


\section{Simple Unital Pseudoalgebras}\lbb{sec.simple}

In this section we classify finitely generated simple unital pseudoalgebra
satisfying certain conditions.  The motivation is to find pseudoalgebras
similar to the most important simple pseudoalgebras, $\cend_n$.

\subsection{General results}
As we are interested in simple unital pseudoalgebras, our objects of study
are necessarily differential (\coref{semisimdiff}).  For
such a pseudoalgebra $R=\dif A$, simplicity easily translates into a
property of the $X^{cop}$-algebra $A$.  Namely, call  $A$ {\em
$X^{cop}$-simple} if it contains no non-zero ideals stable under
the action of $X^{cop}$.  Then we have

\begin{lemma}\lbb{xcopsimple} $\dif A$ is simple if and only if $A$ is
$X^{cop}$-simple.
\end{lemma}

\begin{proof} Let $I$ be an $X^{cop}$-stable ideal of $A$.  Put $\ti
I$ to be the $H$-submodule of $\dif A$ generated by $\{\ti a|\,a\in
I\}$.  Then by (\ref{eq.cend}) and \coref{freeunitdiff}, $\ti I$ is an
ideal of $\dif A$.

Conversely, let $J$ be a non-zero ideal of $\dif A$.  For $b=\sum_K
\d^K\wti{a_K}\in J$, pick $L$ maximal among indices such that $a_L\neq 0$
and consider the product $b_{t^L}e$. It equals $(-1)^L \wti{a_L}$ (see
the proof of \thref{unitisdiff}).  By induction we obtain that all
$\wti{a_K}$ belong to $J$ (compare this to the proof of
\coref{unitsubmod}).  As $\ti a_1\ti b=\wti{ab}$, this implies that $\bar
J=\{a|\,\ti a\in J\}$ is a non-zero ideal of $A$. As $e_x\ti
a=\wti{x(a)}\in J$, this ideal is $X^{cop}$-stable.
\end{proof}

In fact, we proved a more general result:

\begin{lemma}\lbb{idealcorrespondence} The lattice of $X^{cop}$-stable
ideals of $A$ is isomorphic to the lattice of ideals of $\dif A$.
\end{lemma}

\begin{proof} Indeed, in notations from the proof of \leref{xcopsimple}, 
it is clear that $\bar{\ti I}=I$ and $\ti{\bar J}=J$ for any
$X^{cop}$-stable ideal $I$ of $A$ and any ideal $J$ of $\dif A$.
\end{proof}

\subsection{Small $X^{cop}$-Algebras}
Let $A$ be an $X^{cop}$-algebra. We introduce a filtration on $A$:
\begin{equation}\lbb{filtration}
\fil^n A=(\fil_nX)^\perp=\left\{a\in A\,|\,\fil_nX(a)=0\right\}.
\end{equation}

Because of (\ref{dualfilt}) and (\ref{eq.xcop}), the filtration
(\ref{filtration}) respects multiplication in $A$:
\begin{equation}\lbb{eq.filtprop}
(\fil^m A)(\fil^n A)\subset \fil^{m+n} A.
\end{equation}

\begin{remark}\lbb{fullfiltration}  Since the annihilator of every
element of $A$ is not empty, $A=\bigcup_n \fil^n A$.
\end{remark}

We say that a nonzero $a\in A$ has {\em degree} $m$ if
$a\in\fil^m A\backslash \fil^{m-1} A$.  

Remark that the action of $t_i$ lowers the degree, i.e. $t_i(\fil^{m+1}
A)\subset\fil^m A$, otherwise $t_i^{m+1} \fil^{m+1} A\neq 0$.  Notice
also the following useful properties:

\begin{lemma}\lbb{degprop} {\rm (i)} If $\deg a=m$, then for any $M$ such
that $|M|=m$, $t^M(a)\in\fil^0 A$.

{\rm (ii)} If $\deg a=m$, there exists $M$ with $|M|=m$ such that
$t^M(a)\neq 0$.

{\rm (iii)} $\deg(ab)\leq\deg(a)+\deg(b)$.
\end{lemma}

\begin{proof} (i), (ii) follow immediately from (\ref{filtration});
whereas (iii) is a reformulation of (\ref{eq.filtprop}).
\end{proof}

By definition, for $H$ (or, rather, for $H^{op}$) this filtration
coincides with the canonical one.  Here every filtration component is
finite-dimensional.

\begin{lemma}\lbb{finitecomponents} Let $A$ be an $X^{cop}$-algebra such
that $\dim\fil^0A<\infty$.  Then every filtration component is
finite-dimensional.
\end{lemma}

\begin{proof} For any $a\in\fil^{m+1}A\backslash\fil^m A$, there exists
$t_i$ such that $t_i(a)\in\fil^m A\backslash\fil^{m-1} A$. Hence,
\begin{equation}\lbb{eq.compunion}
\fil^{m+1} A=\bigcup_{i=1}^n t_i^{-1}(\fil^m A).
\end{equation}

Notice also that filtration (\ref{filtration}) can be introduced for any
$X$-module $M$ as long as all its elements have non-zero annihilators.  
In particular, this is true for submodules $M_i$ of $A$ defined by
$M_i=\left\{a\in A\,|\,t_ia=0\right\}$.

Consider the particular case of $n=1$ (i.e. $X=\CC[[t]]$). To demonstrate
the statement of the lemma, we induct on $m$.  For any two elements
$a,b\in t^{-1}(\fil^m A)=\fil^{m+1} A$ such that $t(a)=t(b)$, we have
$a-b\in\fil^0 A$.  Thus $\dim\fil^{m+1} A\leq (\dim \fil^m A)(\dim \fil^0
A)\leq (\dim \fil^0 A)^{m+1}$ and $\fil^m A$ is finite-dimensional for all
$m$.

Now we induct on $n$.  Let $a,b\in t_i^{-1}(c)$ for some $i$ and
$c\in\fil^m A$.  Then $a-b\in\fil^m M_i$. Hence (\ref{eq.compunion})
implies:
\begin{align*}
\dim\fil^{m+1} A&\leq\textsum_{i=1}^n \dim t_i^{-1}(\fil^m A)\\
&\leq \textsum_{i=1}^n (\dim\fil^m A)(\dim\fil^m M_i).
\end{align*}
Since by induction $\dim\fil^m M_i<\infty$, we see that $\fil^{m+1} A$ is
finite-dimensional.
\end{proof}

\begin{remark}\lbb{conformalsmall} One of the results in \cite{Re1} can be
reinterpreted as follows: when $H=\CC[\d]$, every finitely generated
$X^{cop}$-simple $X^{cop}$-algebra of $\gk$ (Gelfand-Kirillov dimension,
see \cite{KL}) not exceeding $1$ has finite filtration components.

However, this is not true in general and, by itself, limiting $\gk A$ does
not imply any similarity of $\dif A$ to $\cend_n$ as we will see
immediately below. Remark, first, that for a current pseudoalgebra $\cur
A$, the filtration is trivial: $\fil^0 A=A$.  Unlike the case of conformal
algebras ($H=\CC[\d]$), for larger $H$ there exist non-finite finitely
generated simple current algebras $\cur A$ such that $\gk A\leq\gk H$,
see e.g. \exref{ex.weyl}.  Moreover, consider the following example:

\begin{example}\lbb{ex.nonsmall}  Let $H=U(\gg)$ where $\gg$ is the
three-dimensional abelian algebra.  Then $X=\CC[[t_1,t_2,t_3]]$ and
$X=X^{cop}$.  Let $A=\CC[x,y]$ with the action of $X$ defined by
$t_1=\d/\d x, t_2=0, t_3=0$.  Then $\fil^0 A=\CC[y]$ and we obtain an
$X^{cop}$-algebra that is smaller than $H$ but has a
non-trivial filtration with infinite filtration components.  One can
replace $A$ with the Weyl algebra $A_1=\CC\langle x,y\,|\,xy-yx=1\rangle$
(see \exref{ex.weyl}) and obtain an $X^{cop}$-simple algebra smaller
than $H$ that has a non-trivial filtration with infinite filtration
components. 
\end{example}

Thus, pseudoalgebras similar to $\cend_n$ should not be described by a
simple combinatorial condition such as a bound on $\gk$; although, some
sort of a growth restriction should be imposed.
\end{remark}

We arrive at the following definition:
\begin{definition}\lbb{small}
An $X^{cop}$-algebra $A$ is called {\em small} if $\dim\fil^n A<\infty$
for all $n$.
\end{definition}

Thus, \leref{finitecomponents} can be reformulated as

\begin{lemma}\lbb{zerocomp} Let $A$ be an $X^{cop}$-algebra.  If $\fil^0
A$ is finite-dimensional, then $A$ is small.
\end{lemma}

\begin{remark}\lbb{} It will follow from \thref{smallclass} that for a
simple finitely generated pseudoalgebra $\dif A$ with a small $A$, $A$ has
Gelfand-Kirillov dimension not exceeding $\gk H=\dim\gg$.  This is a
generalization of the converse of the first statement in
\reref{conformalsmall}
\end{remark}

\subsection{Digression: Simplicity Conditions for Small
$X^{cop}$-Algebras}\lbb{digression}
In the next subsection we will show at first that under certain conditions
$A$ can be encoded by an associative algebra ($\fil^0 A$) and a certain
Lie algebra acting on it.  These conditions will be automatically
satisfied when $A$ is small and $X^{cop}$-simple and $\fil^0 A$ is simple.
These two statements about simplicity of $A$ or $\fil^0 A$ are closely 
related.

Indeed, let $A$ be a small $X^{cop}$-algebra. Clearly, if $J$ is a
non-zero proper $X^{cop}$-stable ideal of $A$, $J\cap\fil^0 A$ is an
non-zero ideal of $\fil^0 A$ by \leref{degprop},(ii).

\begin{example}\lbb{innerderiv} Let $A=\End_2(\CC)$ with an inner
derivation $\delta=\ad
\bigl( \begin{smallmatrix}
1&1\\0&1
\end{smallmatrix} \bigr)$.
Then $\fil^0 A=\CC+
\CC \bigl( \begin{smallmatrix}
0&1\\0&0
\end{smallmatrix} \bigr)$.
Therefore, $A$ nay be simple while $\fil^0 A$ is not.
\end{example}

This does not happen when $\fil^0 A$ semisimple. In this case, there exist
non-zero idempotents $a,b\in\fil^0 A$ such that $a(\fil^0 A)b=0$.  If
$aAb=0$, then $AaA$ is a proper $X^{cop}$-stable ideal of $A$, and $A$ is
not $X^{cop}$-simple. However, if there exists $c\in A$ such that $acb\neq
0$,then by applying a suitable element of $X$, we may assume that
$acb\in\fil^0 A$.  But as $acb=a(acb)b=0$, we obtain a contradiction.
 
The above example suggests that the case of an inner derivation is
different from the general one.  This is also true on the pseudoalgebra
level, where the action of an inner derivation can be changed to a trivial
one by changing the pseudoidentity \cite{Re1}.  This leads us to
conjecture that whenever the action of $X^{cop}$ is external in some sense
$A$ is $X^{cop}$-simple if and only if $\fil^0 A$ is simple.  When $X$ is
cocommutative, ``external'' should be understood in the sense of
\cite{Kh}: $X\not\hookrightarrow\Der A_f$, where $A_f$ is the Martindale
quotient of $A$.

\subsection{Simple Small $X^{cop}$-Algebras} 
For the rest of this subsection, $A$ will always stand for a simple small
$X^{cop}$-algebra with $\fil^0 A$ simple (although simplicity is not 
always necessary for the statements below to hold).

Assume now that $A$ satisfies the following technical condition:

\begin{condition}\lbb{cond} As an $\fil^0 A$-module, $\fil^1 A$ is 
generated by $1$ and elements $b_i$, $1\leq i\leq r$, where $r\leq n$, 
such that $t_j(b_i)=\delta_{ij}$.
\end{condition}

\begin{remark}\lbb{condrem} If $\fil^1 A$ is generated over $\fil^0 A$ by
elements $b_i$, $i\in\mathfrak{I}$, where
$\mathfrak{I}\subset\{1,\dots,n\}$ such that $t_j(b_i)=\delta_{ij}$, we
can always renumerate $t_i$'s, so that \condref{cond} holds.

Moreover, we will show below (\leref{subalgofg}) that in this case $A$ is
an an algebra over $\CC[[t_1,\dots,t_r]]^{cop}$, i.e. that
$\CC[[t_1,\dots,t_r]]$ is closed under the action of $\De$.
\end{remark}

We will show in \thref{smallclass} that simplicity of $A$ implies the
above condition (if, of course, $A\neq\fil^0 A$).  The proof is simple but
lengthy, hence we delay it and turn to demonstrating the consequences of
\condref{cond}.  If it holds the structure of $A$ is remarkably nice.  
Namely,

\begin{lemma}\lbb{firstcompgen} If \condref{cond} holds, $A$ is generated
by $\fil^1 A$.  In particular, $A$ is finitely generated.
\end{lemma}

\begin{lemma}\lbb{firstcomplie} If \condref{cond} holds, $\Span(\fil^0
A,b_i)$ is a Lie subalgebra $\bb$ of $A^{(-)}$.  
\end{lemma}

\begin{lemma}\lbb{firstcompactsonzerocomp} If \condref{cond} holds,
$[\bb,\fil^0 A]\subset\fil^0 A$.
\end{lemma}

The proofs of the last two statements are immediate and the proof of the
first comes down to solving a system of linear differential equations:

\begin{proof}[Proof of \leref{firstcompactsonzerocomp}] For any 
$a\in\fil^0 A$ and any $j$, $t_j([b_i,a])=[t_j(b_i),a]=0$.
\end{proof}

We can go further and provide a complete description of $\bb$.  By
\leref{firstcompactsonzerocomp}, $\ad b_i$ is a derivative of $\fil^0 A$
which is a finite-dimensional simple algebra.  Hence, it is inner, i.e.
$\ad b_i=\ad c_i$ \cite{Ja}.  We can replace $b_i$ with $b_i-c_i$, then
the span of $b_i$'s will act trivially on $\fil^0 A$.  It follows that for
any $b_i,b_j$, we have $[b_i,b_j]\in\Span_\CC(b_i)+Z(\fil^0 A)$.  More
explicitly,

\begin{lemma}\lbb{subliefirstcomp} If \condref{cond} holds, we can choose
$b_i$'s, so that for any $1\leq i,j\leq r$, 
\begin{equation}\lbb{eq.subliefirstcomp}
[b_i,b_j]=\sum_{k=1}^r c_{ij}^k b_k+a_{ij},\quad \text{where } 
c_{ij}^k, a_{ij}\in\CC,
\end{equation}
and $[b_i,\fil^0 A]=0$ for all $i$.
\end{lemma}

We can finally turn to \leref{firstcomplie}.

\begin{proof}[Proof of \leref{firstcomplie}] For any $a\in\fil^1 A$ and 
any $j$, by (\ref{coprodform})  
$t_j([a,b_i])=[t_j(a),b_i]+[a,t_j(b_i)]\,\mod\fil^0
A$.  Hence, by \leref{firstcompactsonzerocomp}, $t_j([a,b_i])\in\fil^0 A$
for all $j$.
\end{proof}

\begin{proof}[Proof of \leref{firstcompgen}] Let $B$ be the subalgebra of 
$A$ generated by $\bb$ as defined in the statement of 
\leref{firstcomplie}.

We will prove the following three statements simultaneously:

(i) For all $m$, $\fil^m A\subset B$;

(ii) For all $m$ and any $k>r$, $t_k\fil^m A\subset\fil^{m-2}A$;

(iii) For any collection $\{c_i\}_{i=1}^r$ of elements from $\fil^{m-1}A$
such that $t_i(c_j)=t_j(c_i)$ for all $i,j$, there exists $c\in\fil^m A$
such that $t_i(c)=c_i$.

Clearly (i) will imply the statement of the Lemma.

Remark first that by \condref{cond}, for $m=1$ (i) and (ii) hold
automatically. As for (iii), for $m=1$ we let $c=\sum_i c_ib_i$.

We first demonstrate (ii): let $a\in\fil^m A$.  For any $j$,
$t_jt_k(a)=t_kt_j(a)\in t_k\fil^{m-1} A\subset\fil^{m-3} A$ and we are
done.

Now assume by induction that (iii) holds for $m-1$.  Let $\{c_i\}$ be a
collection of elements from $\fil^m A$. Here and below we will always take
$i\leq r$. By additivity of action of $X$ we may assume $\deg c_i=m$ for
all $i$.  When $X$ is cocommutative, $t_i$'s act simply as derivatives,
thus (iii) comes down to solving a system of linear differential
equations.  Moreover, by induction, first we can pass to $\gr A$, i.e.,
solve the system modulo $\fil^{m-1} A$.  By \leref{subliefirstcomp} it is
equivalent to assuming that $b_i$'s commute with each other and elements
of $\fil^0 A$.  Then the result is classical.

By (\ref{coprodform}), if we consider the action of $t_i$'s modulo
$\fil^{m-1} A$, there is no difference between the general and the
cocommutative case (i.e. the solution for $c$ obtained above is valid
modulo $\fil^{m-1} A$). Hence, there exists an element $c'\in\fil^{m+1} A$
such that $t_i(c')=c_i+d_i$ where $d_i\in\fil^{m-1} A$.  Let $d$ be an
element from $\fil^{m} A$ such that $t_i(d_i)=d$.  Then for $c=c'-d$,
$t_i(c)=c_i$.

We turn to (i).  Assume by induction that $\fil^m\subset B$. For an
arbitrary $a$ of degree $m+1$, put $t_i(a)=a_i\in\fil^m A$.  The
collection $\{a_i\}$ satisfies the conditions of (iii); therefore, we
can produce an element $c\in\fil^{m+1} A$ such that $t_i(c)=a_i$.  By
construction, $c\in B$ as $b_i\in B$ and $a_i,d\in\fil^m A$.  Thus, (ii)
implies that $t_j(a-c)\in\fil^{m-1}A$ for all $j$ and $a-c\in\fil^m A$.  
Therefore, $a\in B$.
\end{proof}

\begin{corollary}\lbb{tkdoesntmatter} If for $a\in A$ and all $i\leq r$, 
$t_i(a)=0$, then $a\in\fil^0 A$.
\end{corollary}

\begin{proof} By \leref{firstcompgen}, $a$ is the sum of monomials of the
type $cb_{i_1}\dots b_{i_m}$, $c\in\fil^0 A$.  Denote $\deg a$ by $m$.
Pick a monomial of degree $m$, say, it ends with $b_j$.  

We may rewrite (\ref{coprodform}) as $\De(t_j)=1\ot t_j+t_j\ot 1+\sum_k
t_k\ot y_{jk}+$~summands that have first terms of degree greater than $1$,
where $y_{jk}\in X$ has no constant terms.  Then it is clear that
monomials of the highest degree in $t_j(a)$ come from monomials in $a$ of
degree $m$.  Hence, $\deg t_i(a)=m-1$.

It follows that if $a$ satisfies the statement of the corollary,
$a\in\fil^1 A$.  \condref{cond} forces $a\in\fil^0 A$.
\end{proof}

We conclude that $A$ is generated by a simple associative algebra $\fil^0
A$ and a Lie subalgebra of $A^{(-)}$ that acts trivially on $\fil^0 A$.
Our goal now is to describe this subalgebra, i.e., to explain the
structural constants $c_{ij}^k$ in (\ref{eq.subliefirstcomp}).

Let $i<j\leq r$.  Clearly, $c_{ij}^k=t_k(b_ib_j-b_jb_i)$ for $k\leq r$.
For $k>r$, we put $c_{ij}^k=0$.  Consider now the structural constants of
$\gg$: $[\d_j,\d_i]=\sum d_{ji}^k \d_k$.  By definition, $d_{ji}^k=\langle
t_k,[\d_j,\d_i]\rangle$.

Since $\d_i\d_j$ is an element of the PBW-basis of $H$, $\langle t_k,
\d_i\d_j\rangle=0$, and we have $d_{ji}^k=\langle t_k,\d_j\d_i\rangle=
\langle {t_k}_{(1)},\d_j\rangle\langle{t_k}_{(2)},\d_i\rangle.$
Therefore, the only summand of $\De(t_k)$ proportional to $t_j\ot t_i$ is
$d_{ji}^kt_j\ot t_i$ (and if $d_{ji}^k=0$, there is no such summand).
Remark also that $\De(t_k)$ has no summand proportional to $t_i\ot t_j$,
otherwise $\langle t_k,\d_i\d_j\rangle\neq 0$.

Thus, comparing expressions for $c_{ij}^k$ and $d_{ji}^k$, we see that
$c_{ij}^k=d_{ji}^k$.

\begin{lemma}\lbb{subalgofg} If \condref{cond} holds,
$\Span(\d_1,\dots,\d_r)$ is a Lie subalgebra of $\gg$.
\end{lemma}

In the same way as above it is not difficult to prove that for $k>r$ in
(\ref{coprodform}), whenever both $t^{K_j}, 
t^{L_j}\in\CC[[t_1,\dots,t_r]]$, $c_j=0$.  Thus, by induction, using the
formula for $\De(t_k)$ stated in the proof of \coref{tkdoesntmatter}, we
have

\begin{lemma}\lbb{tkis0} For $k>r$, $t_k$ acts as $0$ on $A$.
\end{lemma}

We can also strengthen the statement of \leref{degprop}(ii):

\begin{corollary}\lbb{degprop2} If $\deg a=m$, there exists a unique $M$
with $|M|=m$ such that $t^M(a)\neq 0$.
\end{corollary}

Now denote subalgebra $\Span(\d_1,\dots,\d_r)$ by $\hh$.  We can pass to
$H^{op}$ and consider its Lie subalgebra also spanned by
$\d_1,\dots,\d_r$; denote it by $\hh^{op}$. Taking into account that
$a_{ij}$'s in (\ref{eq.subliefirstcomp}) need not be $0$, we obtain

\begin{lemma}\lbb{subalgofgina} If \condref{cond} holds, $b_i$'s generate
a Lie subalgebra of $A$ isomorphic either to $\hh^{op}$, $\hh\subset\gg$,
or its non-trivial $1$-dimensional abelian extension.
\end{lemma}

\begin{theorem}\lbb{smallclass} Let $A$ be an $X^{cop}$-simple small
$X^{cop}$-algebra such that $\fil^0 A$ is simple.  Then $A$ is isomorphic 
to either of

$\ \ (\star)$ $\End_n(\CC)$ with a trivial $X^{cop}$-action;

$\ (\star\star)$ $\End_n(\CC)\ot U(\hh^{op})$, where $\hh$ is a Lie
subalgebra of $\gg$, and the action of $X^{cop}$ is determined by the
action of $U(\hh)^*$;

$(\star\star\star)$ $\End_n(\CC)\ot\big(U(\widehat{\hh^{op}})/(1-{\tt
c})\big)$, where $\widehat{\hh^{op}}$ is a $1$-dimensional abelian
extension $1\to\CC{\tt c}\to\hat{\hh^{op}}\to\hh^{op}\to 1$ of $\hh^{op}$
for a Lie subalgebra $\hh$ of $\gg$, and the action of $X^{cop}$ is
determined by the action of $U(\hh)^*$.

Moreover, in the last two cases $A$ is a simple
$(U(\hh)^*)^{cop}$-algebra.
\end{theorem}

\begin{proof} $\fil^0 A=\End_n(\CC)$.  If $A\neq\fil^0 A$, assume that 
\condref{cond} holds.  Then \leref{subalgofgina} implies that $\fil^1 A$ 
is isomorphic to either $\End_n(\CC)\ot\hh$ or 
$\big(\End_n(\CC)\ot\hat\hh\big)/(1-{\tt c})$, and
by \leref{firstcompgen}, $\fil^1 A$ generates all of $A$.

Therefore, there exists a natural surjective map $\phi$ of either the
algebra of type $(\star\star)$ or $(\star\star\star)$ onto $A$, which is
an isomorphism on the first filtration component.  Notice that these
algebras are small and satisfy \condref{cond}.  Notice also that $\phi$
commutes with the action of $X^{cop}$, in particular, it preserves the
filtration (\ref{filtration}).  To prove injectivity, let $a$ be the
element of least degree such that $\phi(a)=0$.  By statement (iii) of the
proof of \leref{firstcompgen}, there exists $c\in A$ such that
$t_i(c)=t_i\phi(a)$, $i\leq r$.  Let $c'$ be a preimage of $c$ under
$\phi$.  Then by \coref{tkdoesntmatter}, $c'-a$ lie in the zero component
and $\phi(a)\in\phi(c')+\fil^0 A$, a contradiction.

The last claim of the Theorem follows from \leref{tkis0}.

It remains to show that \condref{cond} is valid for simple small
$X^{cop}$-algebras such that $\fil^1 A\neq\emptyset$ and $\fil^0 A$ is 
simple.

Remark that we can change the basis of $\gg$, hence, the generating set 
of $X$.  Thus we will abandon the notation $t_1,\dots,t_n$ that stands
for a fixed generating set and will work with elements of $X$ of degree 
$1$.  Let $T_0=\left\{t\in X\,|\,\deg t=1, t(\fil^1 A)=0\right\}$.

Clearly, for any $t\not\in T_0$, there exists $b\in\fil^1 A, t(b)\neq 0$.
Moreover, if $t\not\in T_0$, there exists $b_t$ such that $t(b_t)=1$.
Indeed, let $b\in\fil^ 1 A$ be such that $t(b)\neq 0$.  Then since $\fil^0
A$ is simple, $s_{1k}, s_{2k}$ such that $t(\sum_k s_{1k}bs_{2k})=\sum_k
s_{1k}t(b)s_{2k}=1$.

In some sense $b_t$, as defined above, is unique.  Put
$N(t)=\{b\,|\,t(b)=0\}$.  Then for any $b$, $b-t(b)b_t\in N(t)$ and $\rk
\fil^1 A/N(t)=1$.

Pick an arbitrary element $t\in X$ of degree $1$, $t\not\in T_0$.  Let
$b_1'$ be such that $t_1(b_1')\neq 0$, and, inductively, $b_j'$ an element
from $\bigcap_{i=1}^{j-1} N(t_i)$ such that there exists $t_j$ for which
$t_j(b_j')\neq 0$.  In this way we obtain sequences $b_1',\dots,b_m'$ and
$t_1,\dots,t_m$.  The process terminates when $\bigcap_{i=1}^{m}
N(t_i)=\fil^0 A$, i.e. when it is annihilated by all $t$'s.  We may assume
that $t_i(b_i')=1$.  Now let $b_m=b_m'$ and, inductively,
$b_i=b_i'-\sum_{j>i} t_j(b_i')b_j$.  In this way we obtain $b_i$'s such
that $t_j(b_i)=\delta_{ij}$ for $i,j\leq m$.

For an arbitrary $b\in\fil^1A\backslash\fil^0 A$, consider the difference
$b-\sum_{i-1}^m t_i(b)b_i$.  As it is killed by all $t_i$'s, by
construction it lies in $\fil^0 A$, hence $1,b_1,\dots,b_m$ is the basis
of $\fil^1 A$ over $\fil^0 A$.  Notice that this is also true when we
consider $\fil^1 A$ as a right $\fil^0 A$-module.

For any $t\not\in T_0$ consider the operator $\sum_{i=1}^m 
t(b_i)t_i$ on $\fil^1 A$.  Using the left and right bases constructed
above, it is easy to see that it acts exactly like $t$.

As in the proof of \leref{firstcompgen}, one can show that
$T_0\fil^m A\subset\fil^{m-2} A$.  Hence, on $\fil^2 A$, $t_jt=t_j(\sum_i
t(b_i)t)$ for any $j$. We can calculate $t_j(\sum_i t(b_i)t_i)(cb_j^2)$,
where $c\in\fil^0 A$, in two ways: either by, applying $t_i$'s first and
then multiplying the results by coefficients $t(b_i)$, or directly by
applying $\De(t)$.  If follows from the discussion immediately preceding
\leref{subalgofg} that $\De(t)(cb_j^2)=t(cb_j)b_j+cb_jt(b_j)$.  By
comparing the results we see that $t(b_j)$ must commute with $c$.
Therefore, $t\in\Span(t_1,\dots,t_m)$ and we are done.
\end{proof}

Clearly, if a simple pseudoalgebra $\dif A$ is finite as an $H$-module,
$A$ is finite, hence it must be of type $(\star)$.  Since $X^{cop}$
acts trivially on $A$, $\dif A$ is necessarily a current algebra.

\begin{corollary}\lbb{smallclassfin} Let $R$ be a simple
differential $H$-pseudoalgebra that is finite as an $H$-module.  Then
$R=\cur\End_n(\CC)$, $n>0$.
\end{corollary}

Now we will describe $H$-pseudoalgebras $\dif A$ when $A$ is
either of the type $(\star\star)$ or $(\star\star\star)$ as defined in
\thref{smallclass}.  So, let $\hh$ be a Lie subalgebra of $\gg$, and
$H'=U(\hh)$.  We can consider $H'$-pseudoalgebra $R'=\dif_{H'}A$.
Recall that in \exref{ex.weyl} we introduced a notation for such
pseudoalgebras: $\cend_n^\phi$ (see also \reref{trivcocycle}).  A
simple comparison of (\ref{eq.cur}) and (\ref{eq.diff}) shows that $\dif_H
A=\cur^H_{H'}R'$.  Therefore, we conclude:

\begin{corollary}\lbb{smallclassforpseudo} Let $R=\dif A$ be a
simple differential $H$-pseudoalgebra such that $A$ is small and $\fil^0 
A$ is simple.  Then either $R=\cur\End_n(\CC)$ or 
$R=\cur^H_{U(\hh)}\cend_n^\phi$ for a Lie subalgebra $\hh$ of $\gg$, 
$\phi\in H^2(\hh)$, and $n>0$.
\end{corollary}

\subsection{Proof of \thref{mainthclass}}
Let $\dif A$ be a simple associative pseudoalgebra.  In the previous
subsection we classified its underlying algebra $A$ satisfying certain 
conditions.  Our goal here is to ``translate'' this condition into one for 
pseudoalgebras.  We need to study the properties of current subalgebras of 
$\dif A$.

\begin{lemma}\lbb{maxcurr}  Let $\dif A$ be a unital associative 
pseudoalgebra with the pseudoidentity $e$.  Then among its current unital 
subalgebras whose pseudoidentity is $e$, $\cur\fil^0 A$ is maximal.
\end{lemma}

\begin{proof} Let $R$ be a unital subalgebra of $\dif A$ with 
pseudoidentity $e$.  Pick $a=\sum_I \d^I\wti{a_I}\in R$. For $J$ maximal 
such that $a_J\neq 0$, $a_{t^J}e\in R$.  By induction, all $\wti{a_I}\in 
R$.

Let $R$ be a current unital subalgebra of $\dif A$ with pseudoidentity $e$
(i.e. $R=\cur B$ and $e=\ti 1$ with regard to the canonical $H$-basis of
$\cur B$).  Let $a=\sum_I \d^I\wti{a_I}$ be an element of the canonical
basis of $\cur B$.  For a non-zero $J$, maximal such that $a_J\neq 0$, 
$a_{t^J}e=0$, hence $a_J=0$, and $a=a_0$.  Since for $I>0$, $e_{t^I}a=0$, 
we conclude that $a_0\in\fil^0 A$ and $a\in\cur\fil^0 A$.
\end{proof}

Therefore, if $\dif A$ is simple and its maximal unital current subalgebra 
with the same pseudoidentity is simple and finite, $A$ is 
small and $X^{cop}$-simple, and $\fil^0 A$ is simple.  

By definition, elements $\ti a$ when $a\in\fil^0 A$ form a (unital)
current subalgebra of $\dif A$.  Hence, by \coref{smallclassfin} and 
\coref{smallclassforpseudo}, the pseudoalgebras satisfying the conditions 
of \thref{mainthclass} are precisely the ones listed there.

It remains to show that pseudoalgebras from that list satisfy the 
conditions of the Theorem.  Simplicity follows from \leref{xcopsimple}.  
The maximal unital current subalgebra is finite simple by \leref{maxcurr}.  
This completes the proof of \thref{mainthclass}.

\subsection{Small Pseudoalgebras}
Finitness of filtration components of $A$ should be translated into 
a finitness condition for $\dif A$.  With this in mind, we propose the 
following

\begin{definition}\lbb{defpseudosmall} A unital differential associative
pseudoalgebra is called {\em small} if all its unital current subalgebras 
are finite as $H$-modules.
\end{definition}

\begin{conjecture}\lbb{smallissmall} A unital pseudoalgebra $\dif A$ is 
small if and only if $A$ is a small $X^{cop}$-simple algebra.
\end{conjecture}

Together with the conjectural statement in subsection \ref{digression},
this will imply a stronger version of \thref{mainthclass}:

\begin{corollary}\lbb{classconj}  A simple unital pseudoalgebra such that 
all its unital current subalgebras are finite over $H$ is either of the 
pseudoalgebras from \thref{mainthclass}.
\end{corollary}

The proof of this conjecture will require understanding what unital 
current subalgebras a differential pseudoalgebra may contain and, more 
generally, the description of the structure of unital subalgebras of 
a differential pseudoalgebra.  The latter is hard to ascertain, as there 
are some counterintuitive examples: for instance, a current conformal 
algebra may contain a non-current unital subalgebra \cite{Re2}.



\begin{thebibliography}{KGB}

\bibitem[BD]{BD} A.~Beilinson, V.~Drinfeld, Chiral algebras,
preprint, 2000, available from 
\textsf{http://zaphod.uchicago.edu/$\tilde{\ }$benzvi/}.

\bibitem[BDK]{BDK} B.~Bakalov, A.~D'Andrea, V.~G.~Kac, Theory of
finite pseudoalgebras, \textit{Adv. Math.} \textbf{162} (2001), pp. 
1--140.

\bibitem[Bl]{Bl} R.~Block, Determination of the Differentiably Simple
Rings with a Minimal Ideal, \textit{Ann. of Math} \textbf{90(3)} (1969),
pp. 433--459.

\bibitem[F]{F} E.~Frenkel, Vertex Algebras and Algebraic Curves,
\textit{Sem. Bourbaki, Exp. 875}, in \textit{Asterisque} \textbf{276}
(2002), pp. 299--341.

\bibitem[G]{G} D.~Gaitsgory, Notes on 2D conformal field theory and string
theory, in \textit{Quantum fields and strings: a course for
mathematicians}, Princeton, NJ, 1996, pp. 1017--1089.

\bibitem[HL]{HL} Y.-Z.~Huang, J.~Lepowsky, On the $\mathcal{D}$-module and 
formal-variable approaches to vertex algebras, in \textit{Topics in 
geometry}, Progr. Nonlinear Differential Equations Appl., \textbf{20}, 
Birkh\"auser, 1996, pp. 175--202.

\bibitem[Ja]{Ja} N.~Jacobson, Abstract derivation and Lie algebras,
\textit{Trans. AMS} \textbf{42(2)} (1937), pp.206--224.

\bibitem[Jo]{Jo} A.~Joseph, \textit{Quantum groups and their primitive
ideals}, Springer-Verlag, Berlin, 1995.

\bibitem[K1]{K1} V.~G.~Kac, \textit{Vertex algebras for beginners},
University Lecture Series, \textbf{10}. AMS, Providence, RI, 1996. Second
edition 1998.

\bibitem[K2]{K2} \bysame, Formal distribution algebras and conformal
algebras, in \textit{Proc. XIIth International Congress of Mathematical
Physics (ICMP '97) (Brisbane)}. Internat. Press, Cambridge, MA, 1999, pp.
80--97. {\texttt{q-alg/9709027}}.

\bibitem[Kh]{Kh} V.~K.~Kharchenko, Constants of derivations of prime 
rings (Russian), \textit{Izv. Akad. Nauk SSSR Ser. Mat.} \textbf{45} 
(1981), pp. 435--461.

\bibitem[KL]{KL} G.~R.~Krause, T.~H.~Lenagan, \textit{Growth of Algebras
and Gelfand-Kirillov Dimension}, second edition. AMS, Providence, RI,
2000.

\bibitem[O]{O} V.~Ostrik, Module categories, weak Hopf algebras, and
modular invariants, preprint, 2001 (\texttt{math.QA/0111139}). 

\bibitem[Re1]{Re1} A.~Retakh, Associative conformal algebras of
linear growth, \textit{J. of Algebra} \textbf{237} (2001), pp. 769--788.

\bibitem[Re2]{Re2} A.~Retakh, Semisimple associative conformal algebras of 
linear growth, in preparation.

\bibitem[Ro]{Ro} M.~Roitman, \textit{Universal enveloping conformal
algebras}, Selecta Math. (N.S.) \textbf{6(3)} (2000), pp. 319--345.

\bibitem[Sw]{Sw} M.~Sweedler, \textit{Hopf algebras}. Math. Lecture Note
Series, W.A.Benjamin Inc., New York, 1969.

\end{thebibliography}
\end{document}